\documentclass[a4paper,11pt]{article}
\usepackage[ansinew]{inputenc}
\usepackage[spanish,english]{babel}
\usepackage{amsmath,amsfonts,amssymb,amsthm}
\usepackage{thmtools}
\usepackage{graphicx}
\usepackage{subfig}
\usepackage[usenames,dvipsnames]{color}
\usepackage{lineno}
\usepackage{enumerate}
\usepackage[colorlinks=true,citecolor=black,linkcolor=black,urlcolor=blue]{hyperref}
\usepackage{amscd,latexsym}
\usepackage{multirow}
\usepackage{cite}

\theoremstyle{plain}
\newtheorem{theorem}{Theorem}
\newtheorem{lemma}{Lemma}

\newtheorem*{invariant*}{Invariant}

\newcommand{\cota}{2\left \lfloor \frac n4 \right \rfloor}
\newcommand{\cotas}{\left \lfloor \frac {2n}5 \right \rfloor}
\newcommand{\cotasn}[1]{\left \lfloor \frac {2(#1)}5 \right \rfloor}

\date{}

\title{Paired and semipaired domination in triangulations
}

\author{M. Claverol\thanks{{\tt merce.claverol@upc.edu}. Universitat Polit\`{e}cnica de Catalunya, Barcelona, Spain.}
\and C. Hernando \thanks{{\tt carmen.hernando@upc.edu}. Universitat Polit\`{e}cnica de Catalunya, Barcelona, Spain.}
\and M. Maureso\thanks{{\tt montserrat.maureso@upc.edu}. Universitat Polit\`{e}cnica de Catalunya, Barcelona, Spain.}
\and M. Mora \thanks{{\tt merce.mora@upc.edu}. Universitat Polit\`{e}cnica de Catalunya, Barcelona, Spain.}
\and J. Tejel\thanks{{\tt jtejel@unizar.es}. IUMA, Universidad de Zaragoza, Zaragoza, Spain.}}

\begin{document}

\maketitle


\begin{abstract}
A dominating set in a graph $G$ is a subset $D$ of vertices such that every vertex not in $D$ is adjacent to at least one vertex in $D$. A dominating set $D$ is paired if the subgraph induced by its vertices has a perfect matching, and semipaired if every vertex in $D$ is paired with exactly one other vertex in $D$ that is within distance 2 from it.
The paired domination number, denoted by $\gamma_{pr}(G)$, is the minimum cardinality of a paired dominating set in $G$, and the semipaired domination number, denoted by $\gamma_{pr2}(G)$, is the minimum cardinality of a semipaired dominating set in $G$.
A near-triangulation is a biconnected planar graph that admits a plane embedding such that all of its faces are triangles except possibly the outer face. In this paper, we show that $\gamma_{pr}(G) \le 2 \left \lfloor \frac{n}{4} \right \rfloor$ for any near-triangulation $G$ of order $n\ge 4$, and that with some exceptions,
$\gamma_{pr2}(G) \le  \left \lfloor \frac{2n}{5} \right \rfloor$ for any near-triangulation $G$ of order $n\ge 5$.
\end{abstract}

\section{Introduction}

Let $G=(V, E)$ be a simple graph. A {\em dominating set} in $G$ is a subset $D\subseteq V$ such that every vertex not in $D$ is adjacent to at least one vertex in $D$.
The {\em domination number} of $G$, denoted by $\gamma (G)$, is defined as the minimum cardinality of a dominating set in $G$.
Dominating sets are of practical interest in several areas and have been widely studied in the literature. Since the mid-1980s, hundreds of papers have been published on graph domination and its variants. 

Among the various domination variants and graph families studied in the literature, this paper focuses on paired and semipaired domination in near-triangulations.
A dominating set $D$ in $G$ is a {\em paired dominating set} ({\em PD-set} for short) if the subgraph induced by the vertices in $D$ has a perfect matching \cite{Haynes98paired}. For every edge $uv$ of this matching, we say that $u$ and $v$ are {\em paired} or {\em partners} in $D$;
if there is no confusion about the paired dominating set to which they belong, we simply say that they are paired or that they are partners. The notion of paired domination was introduced by Haynes and Slater~\cite{Haynes95paired, Haynes98paired}, as a model for the problem of assigning security guards that can protect each other.

Semipaired domination is defined as a relaxed version of paired domination~\cite{Haynes2018semipaired}. Given a simple graph $G =(V, E)$, a dominating set $D$ is a {\em semipaired dominating set} ({\em semi-PD-set} for short) if every vertex in $D$ is matched with exactly one other vertex in $D$ that is within distance 2 from it.
Thus, the vertices in $D$ can be partitioned into pairs such that the distance between the two vertices in each pair is at most 2.
For a pair $(u, v)$, it is said that $u$ and $v$ are {\em semipaired} in $D$ or that they are {\em semipartners} in $D$.

The {\em paired domination number}, denoted by $\gamma_{pr}(G)$, is the minimum cardinality of a PD-set in $G$. Similarly, the {\em semipaired domination number}, denoted by $\gamma_{pr2}(G)$, is the minimum cardinality of a semi-PD-set in $G$. Notice that since the vertices are 
paired or semipaired, the cardinality of a PD-set or a semi-PD-set is always even. In addition, since every PD-set is a semi-PD-set and every semi-PD-set is a dominating set, the following inequalities hold

$$\gamma(G)\le \gamma_{pr2}(G)\le \gamma_{pr}(G)$$
\noindent for every graph $G$.

It is well-known that given a graph $G$ and an integer $k$, the problem of deciding if there is a PD-set in $G$ of cardinality at most $k$ is NP-complete~\cite{Haynes98paired}, even for bipartite graphs~\cite{Chen2010labelling}, split graphs~\cite{Chen2010labelling}, planar graphs~\cite{Tripathi2022complexity} or planar cubic graphs~\cite{Bakal26}. Thus, most of the papers, such as~\cite{Blidia2006characterizations, Chen2009linear, Desormeaux2014, Huang2013paired, Chen2010labelling, Lappas2013n,  lin2020paired, 
Tripathi2022complexity, Bujtas26, Gorz22, Gray25} and the references therein, focus on studying the complexity of computing the paired domination number in several restricted graph classes.

For the semipaired domination, the situation is similar. Given a graph $G$ and an integer $k$, Henning et al.~\cite{Henning20} proved that deciding if there is a semi-PD-set in $G$ of cardinality at most $k$ is NP-complete, even for bipartite or chordal graphs~\cite{Henning20} or planar graphs~\cite{Tripathi25}. In~\cite{Haynes2018totalsemipaired, Haynes2018semipaired, Haynes2019largesemipaired, Henning18claw, Henning19, Henning20, Henning21, Zhuang20, Haynes23, Henning26}, the reader can find recent results on the complexity of computing the semipaired domination number in some particular families of graphs.

One of these families is the well-studied family of maximal outerplanar graphs. In turn,  these graphs are particular cases of near-triangulations\footnote{In the literature, near-triangulations are also called triangulated discs.}. A {\em near-triangulation} is a biconnected planar graph that has a plane embedding such that all of its faces are triangles except possibly the outer face. When the outer face is also a triangle, a near-triangulation is a {\em triangulation} or maximal planar graph.
A {\em maximal outerplanar graph}, MOP graph for short, is a near-triangulation that admits an embedding in which all its vertices lie on the boundary of the outer face.
For a MOP graph $G$ of order $n\ge 4$, Canales et al.~\cite{Canales18} proved that $\gamma_{pr}(G) \le 2 \left \lfloor \frac{n}{4} \right \rfloor$, and for a MOP graph $G$ of order $n\ge 5$, Henning and Kaemawichanurat~\cite{Henning19} showed that $\gamma_{pr2}(G) \le  \left \lfloor \frac{2}{5}n \right \rfloor$, except for a special family $\mathcal{F}$ of MOP graphs of order 9. Both bounds are tight.



In this paper, we extend the results proved in~\cite{Canales18, Henning19} to the family of near-triangulations and show that $\gamma_{pr}(G) \le 2 \left \lfloor \frac{n}{4} \right \rfloor$ for any near-triangulation $G$ of order $n\ge 4$, and that $\gamma_{pr2}(G) \le \left \lfloor \frac{2n}{5} \right \rfloor$ for any near-triangulation $G$ of order $n\ge 5$, except for the special family $\mathcal{F}$ of MOP graphs of order 9 described in~\cite{Henning19}. To our knowledge, these are the first bounds on these two parameters in arbitrary near-triangulations.

The proofs of these two results strongly rely on the techniques given in~\cite{Claverol21} to show that $\gamma_t(G) \le \left \lfloor \frac {2n}{5} \right \rfloor$ for near-triangulations, where $\gamma_t(G)$ is the total domination number of $G$. Given a simple graph $G=(V,E)$, a subset $D\subseteq V$ such that every vertex in $V$ 
is adjacent to a vertex in $D$ is called a {\em total dominating set} ({\em TD-set} for short) in $G$. The {\em total domination number}, denoted by $\gamma _t (G)$, is the minimum cardinality of a total dominating set in $G$.

The paper is organized as follows. Section~\ref{sec:known} is devoted to giving some terminology and basic results on near-triangulations, and to revising the techniques used in~\cite{Claverol21} that we need to prove our main results. In Sections~\ref{sec:paired} and~\ref{sec:semipaired}, we use these techniques to bound the paired and semipaired domination numbers in near-triangulations, respectively. 

\section{Preliminaries}\label{sec:known}

In~\cite{Claverol21}, the authors prove that $\gamma_t(T) \le \left \lfloor \frac {2n}{5} \right \rfloor$ for any near-triangulation $T$, with some exceptions. Their proof is based on a double induction and what the authors call reducible and irreducible near-triangulations, and terminal polygons in irreducible near-triangulations. We will see in the forthcoming sections that
these techniques can also be used to bound the paired and semipaired domination numbers in near-triangulations. Following the definitions and notation in~\cite{Claverol21}, in this section we review these techniques and give some related results.

\begin{figure}[ht!]
	\centering
	\subfloat[]{
		\includegraphics[scale=0.5,page=1]{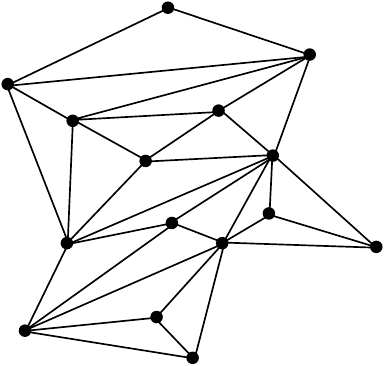}
		\label{fig:NearTriangulationa}
	}~~~~
	\subfloat[]{
		\includegraphics[scale=0.5,page=2]{img-P-SemiP.pdf}
		\label{fig:NearTriangulationb}
    }~~~~
    	\subfloat[]{
		\includegraphics[scale=0.5,page=3]{img-P-SemiP.pdf}
		\label{fig:NearTriangulationc}
    }
	\caption{
		(a) A near-triangulation.
		(b) In bold, the subgraph $T[C]$ induced by the boundary vertices. 
        (c) For the triangular face $uvw$, the vertex $w$ is interior, so the near-triangulation is reducible and the removal of $uv$ gives a near-triangulation.
	}
	\label{fig:NearTriangulations}
\end{figure}

As in~\cite{Claverol21}, we will always refer to a near-triangulation $T=(V, E)$ as a graph that has been drawn in the plane using non-crossing straight-line segments, such that all of its faces are triangles, except possibly the outer face (see Figure~\ref{fig:NearTriangulationa}).
We assume that the boundary of the outer face is a cycle $C=(u_1,u_2,\ldots ,u_h,u_1)$, $h\ge 3$, with its vertices drawn in clockwise order. The vertices and edges of $C$ are \emph{boundary} vertices and edges; a vertex not in $C$ is an \emph{interior} vertex, and an edge connecting two non-consecutive vertices of $C$ is a \emph{diagonal}. Notice that if $h=3$, then $T$ is a triangulation, and if $h=|V|$, then $T$ is a MOP graph.

For any vertex $u\in V$ of a near-triangulation $T=(V,E)$, its {\em open neighborhood} is the set $N_T(u)=\{w\in V: uw\in E\}$ and its {\em closed neighborhood} is $N_T[u]=N_T(u)\cup \{u\}$. We will simply write $N(u)$ or $N[u]$ if the near-triangulation $T$ is clear from context.
The subgraph of $T$ induced by a subset $S$ of vertices in $V$ is denoted by $T[S]$, and
the subgraph obtained from $T$ by removing the vertices $\{v_1, \ldots , v_k\}$ and all their incident edges is denoted by $T-\{v_1, \ldots , v_k\}$. For instance, the bold edges in Figure~\ref{fig:NearTriangulationb} show $T[C]$, the subgraph induced by the boundary vertices.

Similarly, the subgraph obtained from $T$ by removing an edge $e$ is denoted by $T-e$, and  
$T / e$ will denote the graph obtained from $T$ by contracting the edge $e=uv$, that is, the simple graph obtained from $T$ by removing $u, v$ and all their incident edges, adding a new vertex $w$ and connecting $w$ to each vertex $z\in N(u)\cup N(v)$.
We will say that an edge $e$ of $T$ is {\em contractible} if the graph $T / e$ is a near-triangulation.

Several lemmas in~\cite{Claverol21} show how to obtain new near-triangulations when removing vertices and/or edges from a near-triangulation.
The following lemma summarizes some of those cases, which are enough for our purposes.

\begin{lemma}[\cite{Claverol21}]\label{lem:RemovingPoints}
	Let $T$ be a near-triangulation of order $n\ge 4$, with no diagonals, and let $C=(u_1,\ldots ,u_h,u_1)$ be its boundary cycle. 
	Then,
	\begin{enumerate}[i)]
		\item $T-\{u_i\} $ is a near-triangulation for every  $i\in\{1,\dots ,h\}$.
		\item If a boundary edge $e_i=u_{i-1}u_i$ is not contractible in $T$, then there exists an interior vertex $v_i$ adjacent to $u_i$ such that $T-\{ u_i,v_i\}$ is a near-triangulation.
	\end{enumerate}
\end{lemma}

The following lemma, also proved in~\cite{Claverol21}, shows that any boundary vertex is incident to at least one contractible edge.

\begin{lemma}[\cite{Claverol21}] 
	Let $T$ be a near-triangulation of order $n\ge 5$, with boundary cycle $C=(u_1,\ldots ,u_h,u_1)$, and let $u_i$ be a vertex in $C$. Then,	
	\begin{itemize}
		\item[i)]  \label{lem:contractibleedge}  If $u_i$ has a neighbor not in $C$, then there exists an interior vertex $v$ such that the edge $u_iv$ is contractible. 
		\item[ii)]  \label{lem:contractibleNew} If all neighbors of $u_i$ are in $C$, then the edges $u_{i-1}u_i$ and $u_iu_{i+1}$ are contractible. 
	\end{itemize}
\end{lemma}

We say that a near-triangulation $T$ is {\em reducible}~\cite{Claverol21} if it has a boundary edge $e=uv$ such that the third vertex $w$ of the bounded triangular face limited by $e$ is an interior vertex (see Figure~\ref{fig:NearTriangulationc}). Clearly, if we remove $e$ from $T$, then $T-e$ is a near-triangulation of the same order with fewer interior vertices. In addition, $\gamma _{pr}(T)\le \gamma _{pr}(T-e)$, since any PD-set in $T-e$ is also a PD-set in $T$,
and  $\gamma _{pr2}(T)\le \gamma _{pr2}(T-e)$, since any semi-PD-set in $T-e$ is also a semi-PD-set in $T$.

By the definition of reducible near-triangulations, MOP graphs are not reducible near-triangulations because they have no interior vertices.
We say that a near-triangulation is {\em irreducible}~\cite{Claverol21} if it is not reducible and has at least one interior vertex.
Thus, a near-triangulation $T$ is reducible, or irreducible, or a MOP graph.
The smallest irreducible near-triangulation consists of 7 vertices: a triangle with an interior vertex, and three additional vertices of degree 2 added to the sides of the triangle~\cite{Claverol21}.
As a consequence, if $T$ is an irreducible near-triangulation of order $n$, then  $n\ge 7$.

\begin{figure}[tb]
	\centering
	\includegraphics[scale=0.44,page=4]{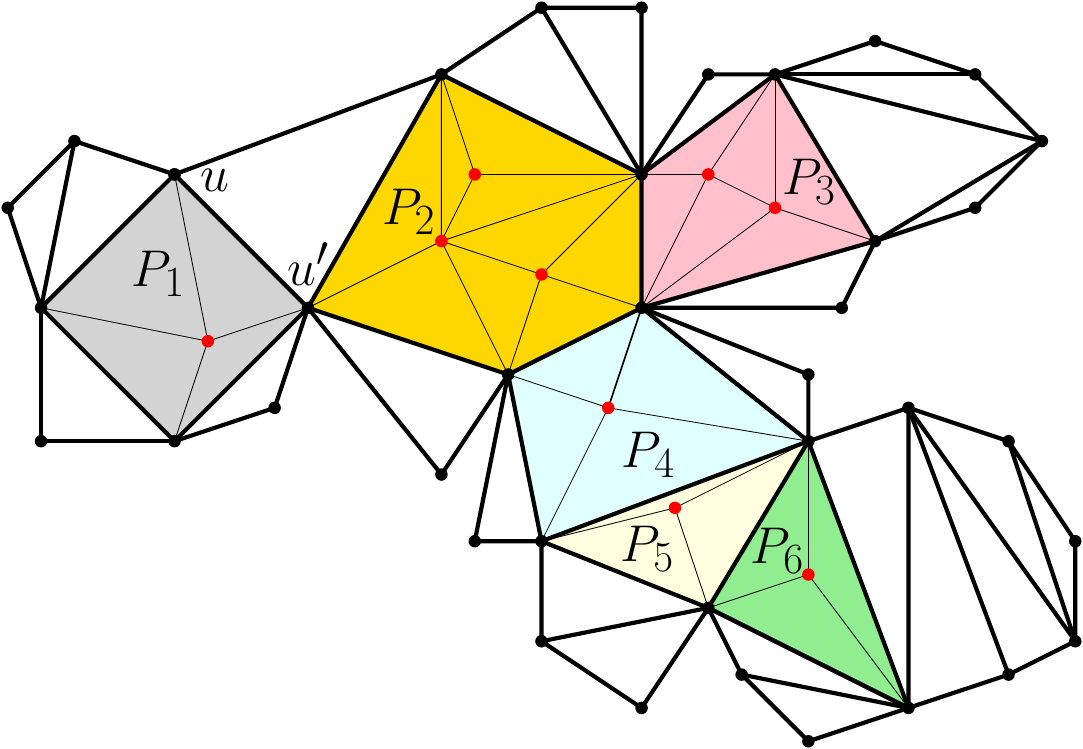}
	\caption{An irreducible near-triangulation $T$. In bold, the subgraph $T[C]$, whose diagonals define a set of adjacent polygons. Six of these polygons are non-empty, $P_1, P_2, P_3,P_4,P_5$ and $P_6$, and three are terminal, $P_1,P_3$ and $P_6$.  The side $(u,u')$ of $P_1$ divides $T$ into two near-triangulations $T_{in}(P_1,(u,u'))$ and $T_{out}(P_1,(u,u'))$, where $T_{in}(P_1,(u,u'))$ is the one containing $P_1$.}\label{fig:Nonreducible}
\end{figure}

Let $T$ be an irreducible near-triangulation $T$ with boundary cycle $C$. The diagonals of the subgraph $T[C]$ divide the interior of $C$ into several regions with disjoint interiors. These regions are simple polygons, whose vertices are boundary vertices of $T$, 
that may or may not contain interior vertices.
For example, the near-triangulation shown in Figure~\ref{fig:Nonreducible} only contains six non-empty polygons, $P_1, P_2, P_3, P_4, P_5$ and $P_6$. The remaining polygons are empty triangles.

Let $P$ be one of these non-empty polygons. Notice that by construction, $P$ has no diagonals of $T$. If $d$ is one of its sides, then $d$ is a diagonal of $T[C]$ and divides $T$ into two non-empty near-triangulations sharing $d$, $T_{in}(P,d)$ and $T_{out}(P,d)$, where $T_{in}(P,d)$ denotes the near-triangulation containing $P$.
It is said that $P$ is {\em terminal}~\cite{Claverol21}, if at most one of the near-triangulations $T_{out}(P,d)$ corresponding to the sides $d$ of $P$
contains interior vertices. Hence, if $P$ is a terminal polygon with $k$ sides, $d_1,\dots,d_k$, then at least $k-1$ of the near-triangulations $T_{out}(P,d_j)$ are MOP graphs with at least three vertices. For instance, the near-triangulation shown in Figure~\ref{fig:Nonreducible} contains three terminal polygons, $P_1, P_3$ and $P_6$.

The following lemma, proved in~\cite{Claverol21}, shows that an irreducible near-triangulation always contains terminal polygons.

\begin{lemma}[\cite{Claverol21}] \label{lem:terminal}
Let $T$ be an irreducible near-triangulation.
Then, $T$ contains at least one terminal polygon.
\end{lemma}

With all these ingredients, we can sketch the proof given in~\cite{Claverol21} that $\gamma_t(T) \le \left \lfloor \frac {2n}{5} \right \rfloor$, for any near-triangulation $T$ (with some exceptions). Let $T$ be a near-triangulation of order $n$ with $m$ interior vertices. The proof that $\gamma_t(T) \le \left \lfloor \frac {2n}{5} \right \rfloor$ follows by induction, assuming the next inductive hypothesis: the bound $\left \lfloor \frac {2n}{5} \right \rfloor$ holds for every near-triangulation $T'$ of order $n'$ with $m'$ interior vertices such that either $n'<n$, or $n=n'$ and $m'<m$ (with some exceptions). If $T$ is reducible, there is an edge $e$ such that $T-e$ is a near-triangulation of order $n$ with $m-1$ interior vertices, 
allowing the inductive hypothesis to be applied to $T-e$. If $T$ is irreducible, then it contains a terminal polygon $P$ with $k\ge 3$ sides. Depending on the orders of the MOP graphs around $P$, several cases are analyzed. 
In all cases, after removing one or several consecutive MOP graphs along with some additional vertices, we can apply the inductive hypothesis to the resulting near-triangulation $T^{\prime }$ and construct specific TD-sets to totally dominate the removed components.
All these sets together define a TD-set in $T$. In addition, a special analysis must be performed when $T'$ is one of the exceptions. As we will see later, the proofs of Theorems~\ref{the:boundpaired} and~\ref{the:boundsemipaired} follow this proof scheme, adapting it to the new scenarios.

We finish this section with some well-known results on MOP graphs. The first, which was already proved 
in~\cite{Shermer91}, shows that a MOP graph can always be divided into two MOP graphs of suitable orders. For the sake of completeness, we include a brief proof, as we add an extra constraint on a boundary edge when partitioning the MOP graph.   

\begin{lemma}\label{lem:diagonales}
	Let $l$ be a positive integer. If $T$ is a MOP graph of order 
    $n\ge 2l$, and $u_iu_{i+1}$ is a boundary edge of $T$, then there exists a diagonal of $T$ that divides $T$ into two MOP graphs, one of them not containing $u_iu_{i+1}$ and of order $h$,  $h\in \{ l+1,\dots, 2l-1\}$.
\end{lemma}
\begin{proof}
	Let $v$ be the third vertex of the triangle of $T$ having the edge $e=u_iu_{i+1}$ as one of its sides. Let $M_1$ and $M_2$ be the MOP graphs defined by the edges $vu_i$ and $vu_{i+1}$, respectively, which do not contain the edge $e$. Note that if $v=u_{i-1}$ (or $v=u_{i+2}$), then $M_1$ ($M_2$) consists of the edge $vu_i$ ($vu_{i+1}$). Let $r_1$ and $r_2$ be the orders of $M_1$ and $M_2$, respectively. We may assume that $r_1\le r_2\le n-1$.
Then $r_1+r_2=n+1\ge 2l+1$, and hence,  $l+1\le r_2\le n-1$. If $r_2\le 2l-1$, the MOP graph $M_2$ satisfies the desired conditions. Otherwise, $r_2\ge 2l$ and we can repeat the reasoning 
for $M_2$ and the boundary edge $vu_{i+1}$ in $M_2$.
Iterating this process yields, in a finite number of steps, a MOP graph of order \(h \in \{l+1, \dots, 2l-1\}\) that does not contain $e$.
\end{proof}

The following lemmas involve TD-sets for MOP graphs of small order.

\begin{lemma}[\cite{Dorfling16, Lemanska17}] \label{lem:pentagon}
Let $T$ be a MOP graph of order 5. For every vertex $u$ of $T$ there exists a TD-set $D$ in $T$ of size 2 containing $u$.
\end{lemma}

\begin{figure}
	\centering
	\includegraphics[scale=0.65,page=6]{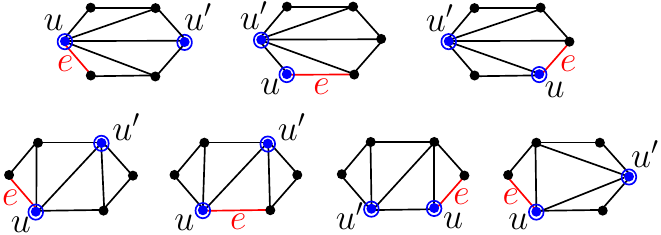}
	\caption{The TD-sets in Lemma~\ref{lem:hexagon}.}
	\label{fig:TDSMOP6}
\end{figure}

\begin{lemma}\label{lem:hexagon}
	Let $T$ be a MOP graph of order 6.
	For every boundary edge $e$ of $T$ there exists a TD-set $D$  in $T$ of size 2 containing {exactly} one of the endpoints of $e$.
\end{lemma}
\begin{proof}
	By symmetry, it is sufficient to analyze the cases of Figure~\ref{fig:TDSMOP6}, 
    which provide TD-sets satisfying the conditions of the statement.
\end{proof}

\begin{lemma}[\cite{Dorfling16, Lemanska17}] \label{lem:heptagon}
Let $T$ be a MOP graph of order 7. There exists a TD-set $D$ in $T$ of size 2.
\end{lemma}

Observe that the two vertices in $D$ of the preceding lemmas must be adjacent, so $D$ is also a PD-set and a semi-PD-set, by pairing or semipairing its vertices.

\section{Paired domination in near-triangulations}\label{sec:paired}

In this section, we prove that $\gamma_{pr}(T) \le 2 \left \lfloor \frac{n}{4} \right \rfloor$ for any near-triangulation $T$ of order $n\ge 4$. 
This result is proved by following the approach in~\cite{Claverol21}, which was sketched in the previous section.
The main differences between the two proofs are the following: (1) now there are fewer cases to analyze, (2) since there are no exceptions, some reasoning can be simplified, and (3) a more careful analysis is required in general to ensure that a perfect matching among the vertices in the dominating set always exists.

We first recall that the same bound $2 \left \lfloor \frac{n}{4} \right \rfloor$ was already proved for MOP graphs by Canales et al.~\cite{Canales18}.

\begin{theorem}[\rm \cite{Canales18}]\label{the:canales}
If $T$ is a MOP graph of order $n\ge 4$, then $\gamma_{pr}(T) \le 2 \left \lfloor \frac{n}{4} \right \rfloor$ and the bound is tight.
\end{theorem}

\begin{theorem}\label{the:boundpaired}
If $T=(V,E)$ is a near-triangulation of order $n\ge 4$, then $\gamma _{pr}(T)\le 2 \left \lfloor \frac{n}{4} \right \rfloor$ and the bound is tight.
\end{theorem}

\begin{proof}
First, notice that the bound is tight, as it has already been established for MOP graphs in~\cite{Canales18}. To prove that the bound holds for any near-triangulation, we proceed by induction on the order and the number of interior vertices of near-triangulations.

Let $T$ be a near-triangulation of order $n\ge 4$, with $m\ge 0$ interior vertices. 
For $m=0$ or $n\le 6$, we claim that the result is true. Indeed, by Theorem~\ref{the:canales}, the upper bound holds for MOP graphs, that is, whenever $m=0$. On the other hand, if $m\not= 0$ and $n\le 6$, then any triangulation is reducible because an irreducible near-triangulation requires at least seven vertices. Thus, since a suitable boundary edge can be removed from a reducible near-triangulation to obtain another near-triangulation with one fewer interior vertex, we can obtain a MOP graph $T'$ of order $n\le 6$, 
by successively removing boundary edges from $T$. Therefore,   $\gamma_{pr}(T)\le \gamma_{pr} (T')\le 2\left \lfloor \frac n4\right \rfloor$, as claimed.

Now, let $T$ be a near-triangulation of order $n$, $n\ge 7$, with $m$  interior vertices, $m>0$.
We want to prove that the result holds for $T$ assuming the following inductive hypothesis: The bound holds for every near-triangulation $T'$ of order $n'$ with $m'$ interior vertices such that either $n'<n$, or $n=n'$ and $m'<m$.

Suppose first that $T$ is reducible. Hence, by removing a suitable boundary edge $e$, we obtain a near-triangulation $T'$ of order $n$ with $m-1$ interior vertices. By the inductive hypothesis,  $\gamma _{pr}(T') \le 2\left \lfloor \frac n4 \right \rfloor$.  Therefore, $\gamma _{pr}(T) \le \gamma _{pr}(T') \le \cota$ and the result holds.

\begin{figure}[htb]
	\centering
	\includegraphics[scale=0.44,page=5]{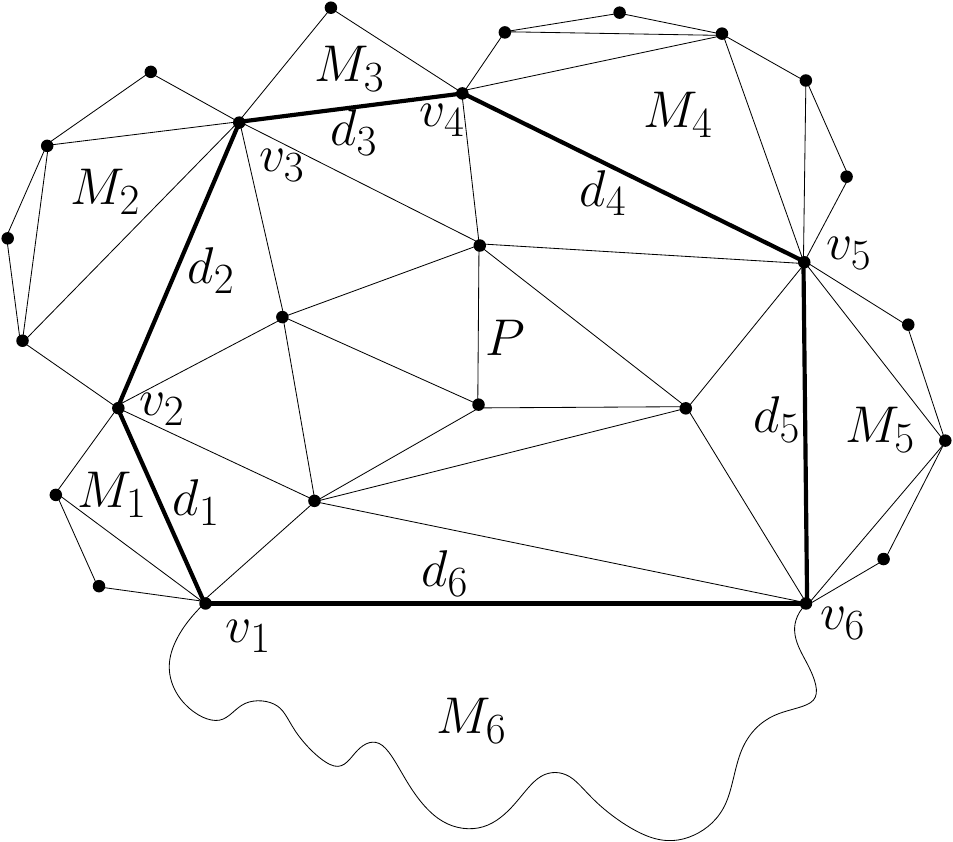}
	\caption{In bold, a terminal $6$-gon $P$ with five MOP graphs $M_1, M_2, M_3, M_4$ and $M_5$ of orders $4,6,3,7$ and $5$, respectively, around it. $M_6$ could be a MOP graph or not.}\label{fig:TerminalPolygon}
\end{figure}

Suppose now that $T$ is irreducible.
By Lemma~\ref{lem:terminal}, $T$ contains a terminal polygon $P$ with $k\ge 3$ sides, $d_1=v_1v_2$, $d_2=v_2v_3$, $\dots $, $d_k=v_kv_1$. We use $T_j$ to denote $T_{in}(P,d_j)$ and $M_j$ to denote $T_{out}(P,d_j)$, for $j=1,\ldots ,k$. Obviously, $|T_j|+|M_j|=n+2$. In addition, since $P$ is a terminal polygon, $M_j$ is a MOP graph for every $j\in \{ 1,\dots ,k\}$, except possibly one of them, say $M_k$. Figure~\ref{fig:TerminalPolygon} illustrates an example of a terminal 6-gon.
Observe that since $T$ is irreducible, we have $|T_j| \ge 6$, for every $j\in \{1,\ldots ,k\}$.
In addition, since $P$ has interior vertices and contains no diagonals, the vertex $w$ in the triangular face $v_jv_{j+1}w$ of $T_j$ must be an interior vertex, so $T_j$ is reducible, and $T_j-d_j$ is a near-triangulation.

We distinguish cases according to the orders of the MOP graphs $M_j$, $j\in \{1,\dots , $ $k-1\}$.
Concretely,
if there is a MOP graph $M_j$ of order 4, 5 and $d_j$ is contractible in $T_j$, 6, 7, or at least 8 (Cases 1-5); or if all MOP graphs have order 3 or 5, with $d_j$ not contractible in $T_j$ when the order is 5 (Cases 6-8). In cases 1-5, it suffices to remove vertices from a single MOP graph to apply the inductive hypothesis, whereas in cases 6-8, we will remove vertices from two consecutive MOP graphs, $M_j$ and $M_{j+1}$.

\vskip 0.2 cm
\noindent
{\bf Case 1:} $|M_j|=4$, for some $j\in \{1,\dots, k-1\}$.
\vskip 0.1 cm

Suppose that there is a MOP graph $M_j$ of order 4 ($M_1$ in Figure~\ref{fig:TerminalPolygon}), whose vertices are $v_j, w_1, w_2$ and $v_{j+1}$ in clockwise order. Hence, one of the vertices, $v_j$ or $v_{j+1}$, has degree 3 in $M_j$, say $v_j$, and $T_j=T-\{w_1,w_2\}$.

Let $v$ be the boundary vertex in $T_j$ (and in $T$) preceding $v_j$; see Figure~\ref{fig:case1-size4}. Let $T'$ be the near-triangulation obtained by removing $d_j$ from $T_j$ and adding two new vertices $x$ and $y$ to $T_j$ such that $x$ is adjacent to $v$ and $v_{j}$, and $y$ is adjacent to $x$ and $v_{j}$ (see Figure~\ref{fig:case1-size4}). Hence, $T'$ is a near-triangulation of order $n$, with $m-1$ interior vertices.
By the inductive hypothesis, there is a PD-set $D'$ in $T'$ of size at most $\cota $.


From $D'$, we can obtain a PD-set $D$ in $T$ of size at most $\cota $ as follows. We analyze several cases, depending on whether $v_{j}, x$, and $y$ belong to $D'$ or not. Observe that if two vertices different from $x$ and $y$ are paired in $D'$, then they can also be paired in a dominating set of $T$. 
Suppose first that $v_{j}$ belongs to $D'$.
\begin{enumerate}
\item[$\bullet$]  If $\{x,y\}\notin D'$, then $D=D'$ is clearly a PD-set in $T$, as $v_j$ dominates $w_1$ and $w_2$; see Figure~\ref{fig:m4a}.

\begin{figure}[!htb]
	\centering
    \subfloat[The near-triangulation $T'$ obtained from $T$.]{
	\includegraphics[scale=0.46,page=7]{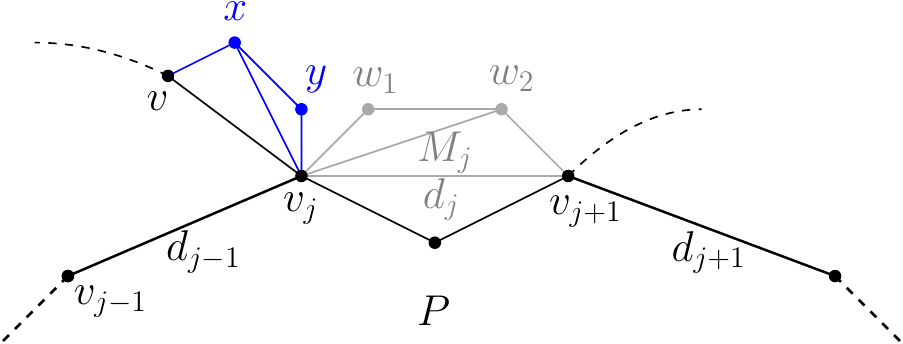}
	\label{fig:case1-size4}
	}\\
	\subfloat[$v_{j}\in D'$ and $\{x,y\}\notin D'$.]{
		\includegraphics[scale=0.46,page=18]{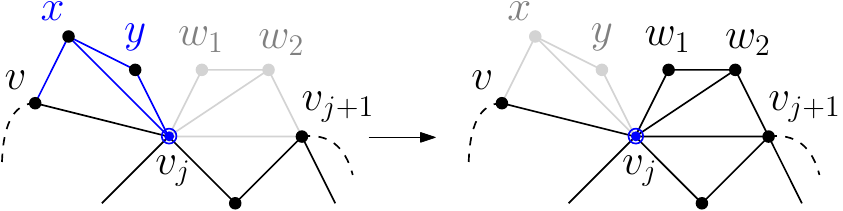}
		\label{fig:m4a}
	}~~~~
	\subfloat[$v_{j}\in D'$, $x\notin D'$ and $y\in D'$.]{
		\includegraphics[scale=0.46,page=19]{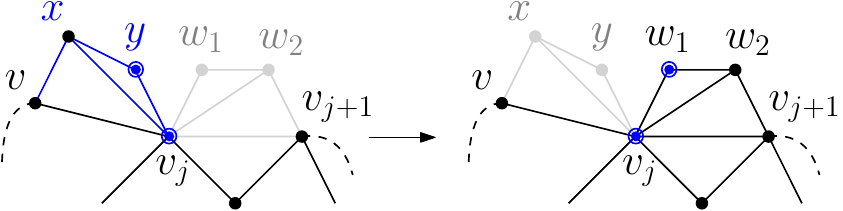}
		\label{fig:m4b}
   }\\
   	\subfloat[$v_j$ and $x$ are paired in $D'$.]{
		\includegraphics[scale=0.46,page=20]{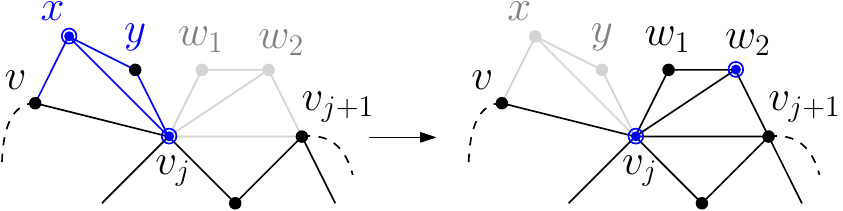}
		\label{fig:m4c}
	}~~~~
	\subfloat[$v_j$ and $y$ are paired in $D'$.]{
		\includegraphics[scale=0.46,page=21]{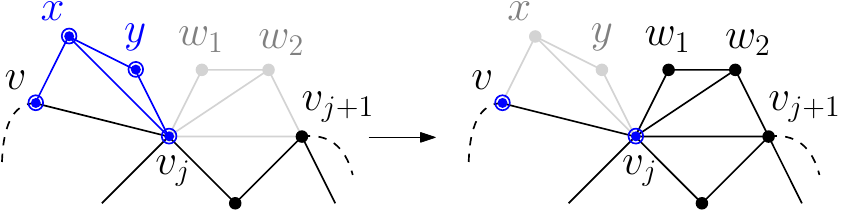}
		\label{fig:m4d}
   }\\
   	\subfloat[$v_{j}\notin D'$, $x\in D'$ and $y\notin D'$.]{
		\includegraphics[scale=0.46,page=22]{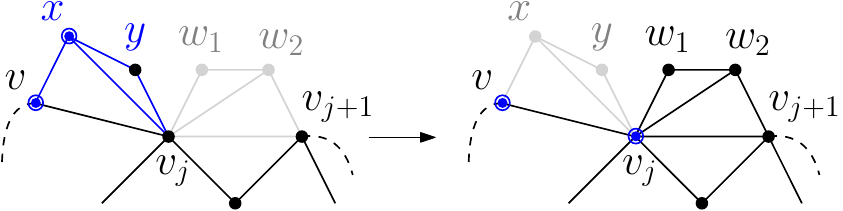}
		\label{fig:m4e}
	}~~~~
	\subfloat[$v_{j}\notin D'$, and $x$ and $y$ are paired in $D'$.]{
		\includegraphics[scale=0.46,page=23]{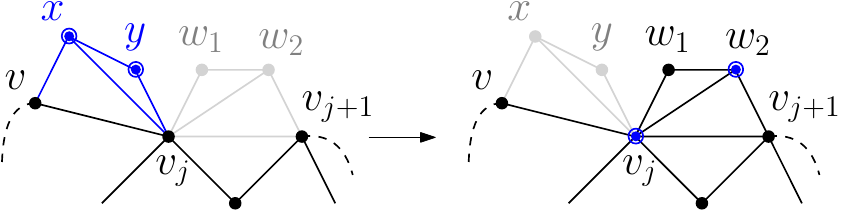}
		\label{fig:m4f}
   }
    \caption{Illustrating Case 1 of Theorem~\ref{the:boundpaired}.}
    \label{fig:m4}
\end{figure}

\item[$\bullet$]  If  $x\notin D'$ and $y\in D'$, then $y$ and $v_{j}$ must be paired in $D'$; see Figure~\ref{fig:m4b}. Hence, $D=(D'\setminus \{y\}) \cup \{w_1\}$ is a PD-set in $T$, where $v_j$ and  $w_1$ are paired in $D$.

\item[$\bullet$]  If $x\in D'$ and $y\notin D'$, then $x$ is paired with either $v_{j}$ or $v$. If  $x$ and  $v_{j}$ are paired in $D'$ (see Figure~\ref{fig:m4c}), then the set $D=(D'\setminus \{x\}) \cup \{w_2\}$ is a PD-set in $T$, with $v_j$ and $w_2$ paired in $D$. Otherwise, $x$ and $v$ are paired in $D'$. We consider whether there exists a neighbor $v'$ of $v$ not belonging to $D'$. If such a vertex exists, then $D=(D'\setminus \{x\}) \cup \{v'\}$ is a PD-set in $T$, where $v$ and $v'$ are paired in $D$. If all the neighbors of $v$ belong to $D'$, then $D=D'\setminus \{x,v\}$ is a PD-set in $T$, as $v, w_1$ and $w_2$ are all neighbors of $v_j$.

\item[$\bullet$] If $x\in D'$ and $y\in D'$, then $y$ is paired with either $x$ or $v_j$ in $D'$. If $x$ and $y$ are paired in $D'$, then $D=D'\setminus \{x,y\}$ is a PD-set in $T$. If $y$ and $v_{j}$ are paired in $D'$ (see Figure~\ref{fig:m4d}), then $x$ and $v$ must be paired in $D'$. Then, $D=D'\setminus\{x,y\}$ is a PD-set in $T$, with $v$ and $v_j$ paired in $D$.
\end{enumerate}

Suppose now that $v_{j}$ does not belong to $D'$.
\begin{enumerate}
\item[$\bullet$]
If $y\notin D'$, then $x$ must belong to $D'$ to dominate $y$, and $x$ must be paired with $v$, as illustrated in Figure~\ref{fig:m4e}. In this case, $D=(D'\setminus \{x\}) \cup \{v_{j}\}$ is PD-set in $T$, where $v$ and $v_j$ are paired in $D$.
\item[$\bullet$] If $y\in D'$, then $x$ must be its partner in $D'$; see Figure~\ref{fig:m4f}. Thus, $D=(D'\setminus \{x,y\}) \cup \{v_{j}, w_2\}$ is a PD-set in $T$, with $v_j$ and $w_2$ paired in $D$.
\end{enumerate}

\vskip 0.2 cm
\noindent
{\bf Case 2:} $|M_j|=5$ and $d_j$ is contractible in $T_j$, for some $j\in \{ 1,\dots, k-1\}$.

\vskip 0.1 cm
Suppose that there is a MOP graph $M_j$ of order 5 ($M_5$ in Figure~\ref{fig:TerminalPolygon}) such that $d_j=v_j v_{j+1}$ is contractible in $T_j$. Let $v_j, w_1, w_2, w_3$ and $v_{j+1}$ be the vertices of $M_j$, clockwise ordered; see Figure~\ref{fig:m5}. Hence, $T_j=T-\{w_1,w_2,w_3\}$ and $|T_j|\ge 6$.
Let $T'=T_j/d_j$. Since $d_j$ is contractible in $T_j$, $T'$ is a near-triangulation of order $n-4$, with $n-4\ge 5$. By the inductive hypothesis, $T'$ has a PD-set $D'$ of size at most $2\left \lfloor \frac{n-4}4\right \rfloor$. In addition, by Lemma~\ref{lem:pentagon}, there is a TD-set $D''$ of size 2 in $M_j$  containing a fixed vertex of $M_j$.
From $D'$ and $D''$, we build a PD-set $D$ in $T$ as follows.

\begin{figure}[!htb]
	\centering
	\subfloat[$D''=\{w_1,w_2\}$.]{
		\includegraphics[scale=0.48,page=15]{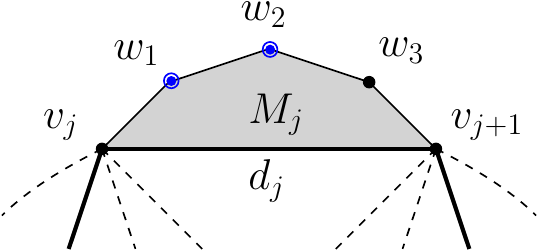}
		\label{fig:m5a}
	}~
	\subfloat[$D''=\{v_j,w_2\}$.]{
		\includegraphics[scale=0.48,page=16]{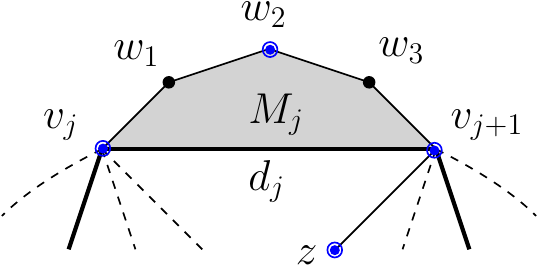}
		\label{fig:m5b}
   }~
   	\subfloat[$v_{j+1}$ belongs to $D''$.]{
		\includegraphics[scale=0.48,page=17]{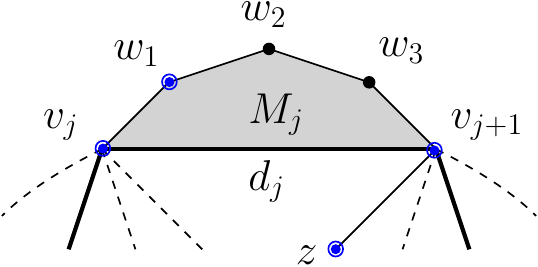}
		\label{fig:m5c}
   }
    \caption{Illustrating Case 2 of Theorem~\ref{the:boundpaired}.}
    \label{fig:m5}
\end{figure}

Let $w$ be the vertex in $T'$ corresponding to the contracted edge $d_j$. Suppose first that $w\notin D'$.  In this case, $D=D'\cup D''$ is clearly a PD-set in $T$ of size at most $2\left \lfloor \frac{n-4}4\right \rfloor+2 = \cota$; see Figure~\ref{fig:m5a} for an illustration. Indeed, the vertices in $D'$ and $D''$ are paired, and a vertex $x$ not in $D'\cup D''$ is dominated either by a vertex in $D''$, if $x$ belongs to $M_j$, or by a vertex in $D'$, if $x\ne v_j, v_{j+1}$ belongs to $T_j$.

Suppose now that $w\in D'$ and let $z$ be its partner in $D'$. Thus, $z$ is adjacent in $T$ to $v_j$ or $v_{j+1}$, say $v_{j+1}$. Let $D''$ be a TD-set in $M_j$ containing $v_j$, which exists by Lemma~\ref{lem:pentagon}. We distinguish two cases.

\begin{enumerate}
\item[$\bullet$] If $v_{j+1}\notin D''$, then $D=(D'\setminus \{w\})\cup \{v_{j+1}\}\cup D''$ is a PD-set in $T$ of size at most $2\left \lfloor \frac{n-4}4\right \rfloor+2 = \cota$, where $v_{j+1}$ and $z$ are paired in $D$. See Figure~\ref{fig:m5b}. Notice that if a vertex $x$ was dominated by $w$ in $T'$, it is now dominated by $v_j$ or $v_{j+1}$ in $T$.
\item[$\bullet$] If $v_{j+1}\in D''$, then $D=(D'\setminus \{w\})\cup \{v_{j+1}, v_j, w_1\}$ is a PD-set in $T$ of size at most $2\left \lfloor \frac{n-4}4\right \rfloor+2 = \cota$, where $v_{j+1}$ is paired with $z$ in $D$, and $v_j$ with $w_1$. See Figure~\ref{fig:m5c}.
\end{enumerate}

\vskip 0.2 cm
\noindent
{\bf Case 3:} $|M_j|=6$, for some $j\in \{1,\dots, k-1\}$.
\vskip 0.1 cm

Suppose that there is a MOP graph $M_j$ of order 6; $M_2$ in Figure~\ref{fig:TerminalPolygon} (see also Figure~\ref{fig:m6}).
Let $v_j, w_1, w_2, w_3,w_4$ and $v_{j+1}$ be the vertices of $M_j$, in clockwise order.
Hence, $T_j=T-\{w_1,w_2,w_3,w_4\}$, with $|T_j|=n-4\ge 6$.
By the inductive hypothesis, there exists a PD-set $D'$ in $T_j$ of size at most $2\left \lfloor \frac{n-4}4 \right \rfloor$.
In addition, by Lemma~\ref{lem:hexagon}, there is a TD-set $D''=\{u,v\}$ of size 2 in $M_j$ that contains exactly one of $v_j$ or $v_{j+1}$, say $u=v_j$. Note that $v\notin T_j$, so either $D'\cap D''=\emptyset$ or $D'\cap D''=\{v_j\}$.

\begin{figure}[!htb]
	\centering
	\subfloat[$D''=\{v_j, w_4\}$ and $v_j \notin D'\cap D''$.]{
		\includegraphics[scale=0.46,page=24]{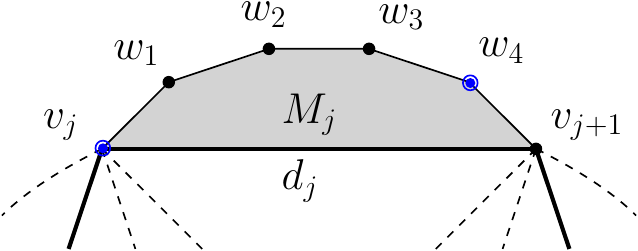}
		\label{fig:m6a}
	}~~~~~~~~~~~~~~~~
	\subfloat[$D''=\{v_j, w_4\}$ and $v_j \in D'\cap D''$.]{
		\includegraphics[scale=0.46,page=25]{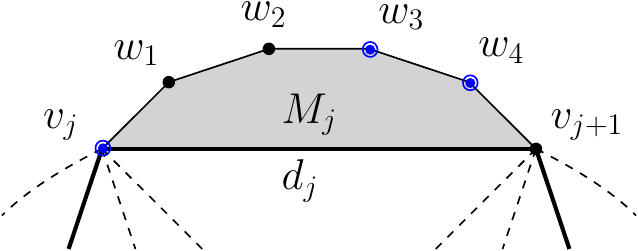}
		\label{fig:m6b}
    }
    \caption{Illustrating Case 3 of Theorem~\ref{the:boundpaired}.}
    \label{fig:m6}
\end{figure}

From $D'$ and $D''$, we build a PD-set $D$ in $T$ as follows. If $D'\cap D''=\emptyset$, then $D=D'\cup D''$ is a PD-set in $T$ of size at most $2\left \lfloor \frac{n-4}4\right \rfloor+2 = \cota$; see Figure~\ref{fig:m6a} for an illustration. If $D'\cap D''=\{v_j\}=\{u\}$, then at least one of the neighbors of $v$ on the boundary of $M_j$ does not belong to $D'\cup D''$ ($w_3$ in Figure~\ref{fig:m6b}). Let $v'$ be such a neighbor. The set $D=D'\cup \{ v,v'\}$ is a PD-set in $T$, where $v$ and $v'$ are paired in $D$, of size at most $2\left \lfloor \frac{n-4}4\right \rfloor+2 = \cota$.

\vskip 2mm

\noindent
{\bf Case 4:} $|M_j|=7$, for some $j\in \{1,\dots, k-1\}$.
\vskip 1mm

Suppose that there is a MOP graph $M_j$ of order 7 ($M_4$ in Figure~\ref{fig:TerminalPolygon}); see also Figure~\ref{fig:m7}.
Since $T_j$ is a near-triangulation of order $n-5\ge 6$, by the inductive hypothesis, $T_j$ has a PD-set $D'$ of size at most $2\left \lfloor \frac{n-5}4 \right \rfloor$.
In addition, by Lemma~\ref{lem:heptagon}, $M_j$ has a TD-set $D''=\{u,v\}$ of size 2.

\begin{figure}[!htb]
	\centering
	\subfloat[$D'\cap D''=\emptyset$.]{
		\includegraphics[scale=0.46,page=26]{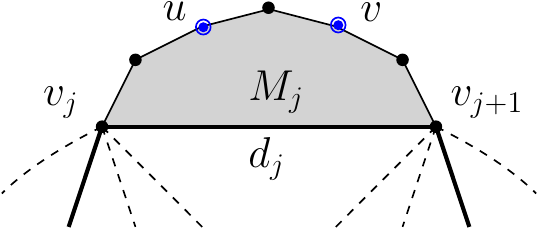}
		\label{fig:m7a}
	}
	\subfloat[$D'\cap D''=\{v_j\} =\{u\}$.]{
		\includegraphics[scale=0.46,page=27]{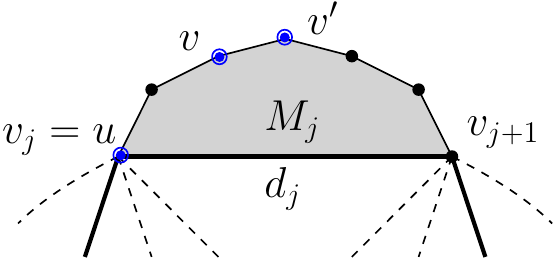}
		\label{fig:m7b}
   }
   	\subfloat[$D'\cap D''=\{u, v\}$.]{
		\includegraphics[scale=0.46,page=28]{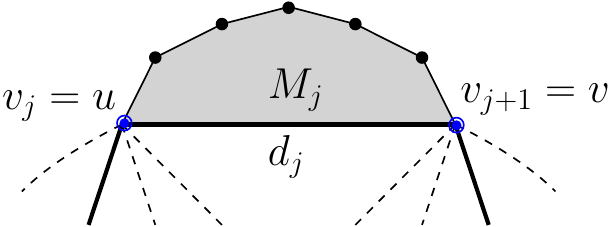}
		\label{fig:m7c}
   }
    \caption{Illustrating Case 4 of Theorem~\ref{the:boundpaired}.}
    \label{fig:m7}
\end{figure}

If $D'\cap D''=\emptyset$, then $D=D'\cup D''$ is a PD-set in $T$; see Figure~\ref{fig:m7a}. If $|D'\cap D''|=1$, say $D'\cap D''=\{v_j\} =\{u\}$, then $D=D'\cup \{ v,v'\}$ is a PD-set in $T$, where $v'$ is one of the neighbors of $v$ on the boundary of $M_j$ that does not belong to $D'\cup D''$; see Figure~\ref{fig:m7b}. Finally, if $D'\cap D''=\{u, v\}$, then $D=D'$ is a PD-set in $T$; see Figure~\ref{fig:m7c}. In any case, $D$ is a PD-set of size at most $2\left \lfloor \frac{n-5}4\right \rfloor+2 \le \cota$.

\vskip 0.2 cm
\noindent
{\bf Case 5:} $|M_j|\ge 8$, for some $j\in \{1,\dots, k-1\}$.
\vskip 0.1 cm

Suppose that there is a MOP graph $M_j$ of order at least 8. Taking $l=4$ in Lemma~\ref{lem:diagonales}, there is a diagonal $d$ in $M_j$ dividing $M_j$ into two MOP graphs, one of them, say $M'$, of order $r\in \{5,6,7\}$ and not containing the edge $d_j$; see Figure~\ref{fig:m9} for an illustration. Thus, 
we can proceed as in Cases 2, 3, and 4 to construct a PD-set in $T$ of size at most $\cota$, taking $M'$ as the corresponding $M_j$. Observe that when removing the vertices of $M'$ (except for the endpoints of $d$) from $T$ to obtain a near-triangulation $T'$, the edge $d$ is a boundary edge in $T'$. Hence, by Lemma~\ref{lem:contractibleNew} $ii)$, $d$ is contractible in $T'$ because at least one of the endpoints of $d$ has all its neighbors on the boundary of $T'$.

\begin{figure}[!htb]
	\centering
	\includegraphics[scale=0.45,page=29]{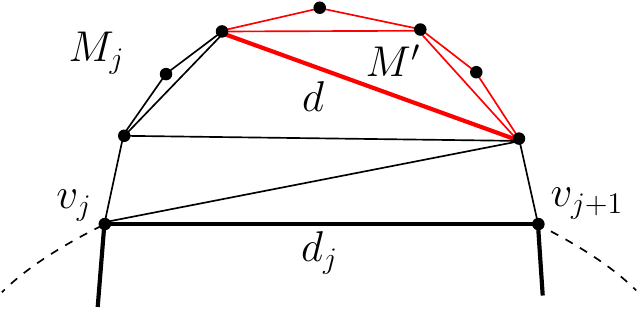}
	\caption{$M_j$ is a MOP graph of order 9. The diagonal $d$ defines a MOP graph $M'$ of order 5.}\label{fig:m9}
\end{figure}

\vskip 0.2 cm
\noindent
{\bf Case 6:} $|M_j|=3$ and $|M_{j+1}|=3$, for some $j\in \{1,\dots,k-2\}$.
\vskip 0.2 cm


\begin{figure}[!htb]
	\centering
	\subfloat[The near-triangulation $T'$ obtained from $T$.]{
		\includegraphics[scale=0.46,page=8]{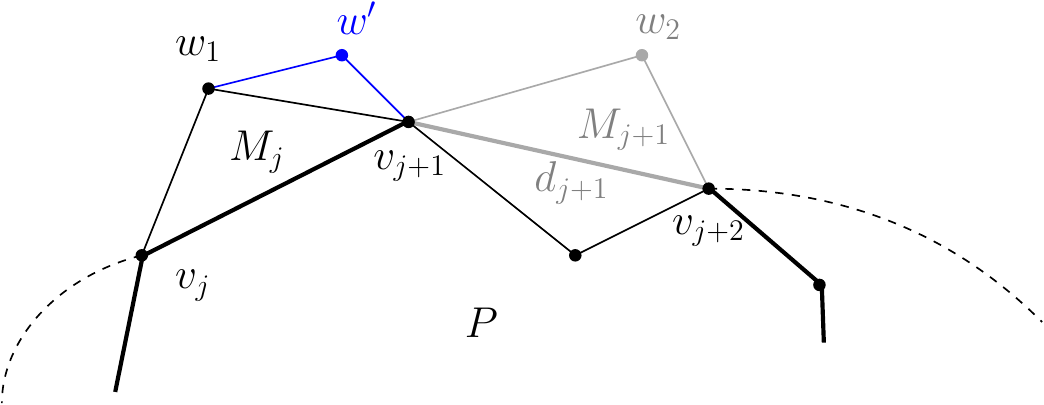}
		\label{fig:case6b}
	}\\
	\subfloat[$w$ and $v_{j+1}$ are paired.]{
		\includegraphics[scale=0.46,page=30]{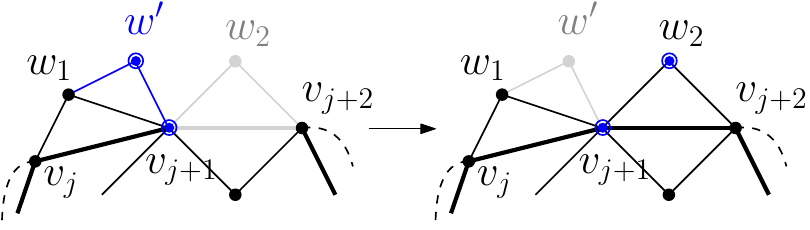}
		\label{fig:case6_1}
   }~~~~
   	\subfloat[$v_{j+1}\in D'$, and $w'$ and $w_1$ are paired.]{
		\includegraphics[scale=0.46,page=31]{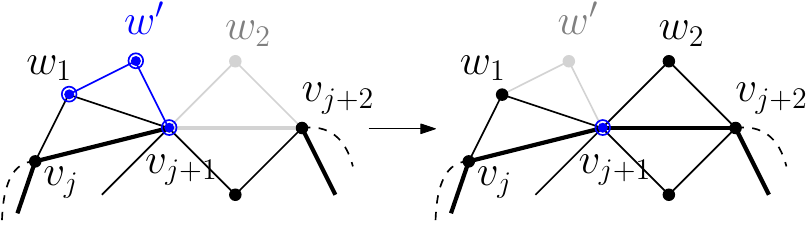}
		\label{fig:case6_2}
        	}\\
	\subfloat[$v_{j+1}\notin D'$, and $w'$ and $w_1$ are paired.]{
		\includegraphics[scale=0.46,page=32]{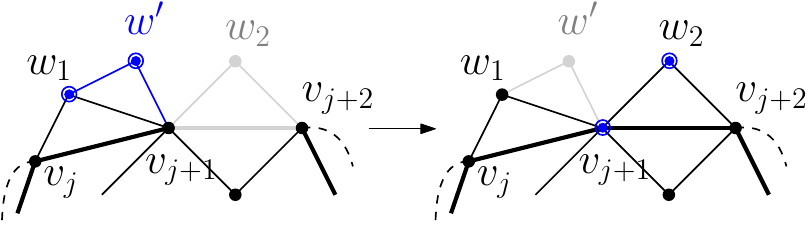}
		\label{fig:case6_3}
   }~~~~
   	\subfloat[$\{w', v_{j+1}\}\notin D'$, and $w_1$ and $v_j$ are paired.]{
		\includegraphics[scale=0.46,page=33]{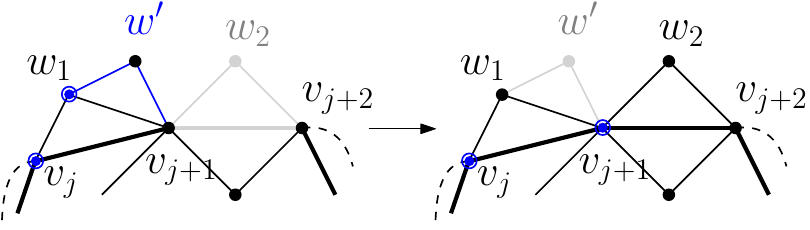}
		\label{fig:case6_4}
   }
    \caption{Illustrating Case 6 of Theorem~\ref{the:boundpaired}.}
    \label{fig:case6}
\end{figure}

Suppose that there exist two consecutive MOP graphs $M_j$ and $M_{j+1}$ of order 3. This case is similar to Case 1. Let $v_j, w_1, v_{j+1}$ and $v_{j+1}, w_2, v_{j+2}$ be the vertices of $M_j$ and $M_{j+1}$, respectively. 
Let $T'$ be the near-triangulation obtained by removing the vertex $w_2$, adding a new vertex $w'$ adjacent to $w_1$ and $v_{j+1}$, and removing the edge $d_{j+1} $. Then, $T'$ is a near-triangulation of order $n$ with $m-1$ interior vertices (see Figure~\ref{fig:case6b}). By the inductive hypothesis, there is a PD-set $D'$ in $T'$ of size at most $\cota$.
From $D'$, we obtain a PD-set $D$ in $T$ of size at most $\cota$. We distinguish two cases, depending on whether $w'$ belongs to $D'$ or not.

Suppose first that $w'$ belongs to $D'$. In this case, $w_1$ or $v_{j+1}$ must belong to $D'$, in order to ensure that $w'$ is paired with $w_1$ or with $v_{j+1}$.
	\begin{enumerate}
		\item[$\bullet$] If $w'$ and $v_{j+1}$ are paired in $D'$, then $D=(D'\setminus\{w'\})\cup \{w_2\}$ is a PD-set in $T$ of size at most $\cota$, where $v_{j+1}$ and $w_2$ are paired; see Figure~\ref{fig:case6_1}.
		\item[$\bullet$] If $w'$ and $w_1$ are paired in $D'$, we distinguish two cases. If $v_{j+1}\in D'$, then $D=D'\setminus\{w',w_1\}$ is a PD-set in $T$ of size at most $\cota-2<\cota$, as all neighbors of $w_1$ are also neighbors of $v_{j+1}$; see Figure~\ref{fig:case6_2}. If $v_{j+1}\notin D'$, then $D=(D'\setminus\{w',w_1\}) \cup \{v_{j+1},w_2\}$ is a PD-set in $T$ of size at most $\cota$; see Figure~\ref{fig:case6_3}.
	\end{enumerate}

Now suppose that $w'$ does not belong to $D'$. If $v_{j+1}\in D'$, then $D=D'$ is a PD-set in $T$ of size at most $\cota$. If $v_{j+1}\notin D'$, then $w_1\in D'$ and $v_j$ must be its partner in $D'$. Hence, $D=(D'\setminus\{w_1\})\cup \{v_{j+1}\}$ is a PD-set in $T$ of size at most $\cota$, where $v_j$ and $v_{j+1}$ are partners; see Figure~\ref{fig:case6_4}.

\vskip 0.2 cm
\noindent
{\bf Case 7:} $|M_j|=5$ with $d_j$ not contractible in $T_j$,
and $|M_{j+1}|=3$,  for some $j\in \{1,\dots,k-2\}$.

\vskip 0.1 cm

\begin{figure}[!httb]
	\centering
	\subfloat[The near-triangulation $T'$ obtained from $T$.]{
		\includegraphics[scale=0.46,page=9]{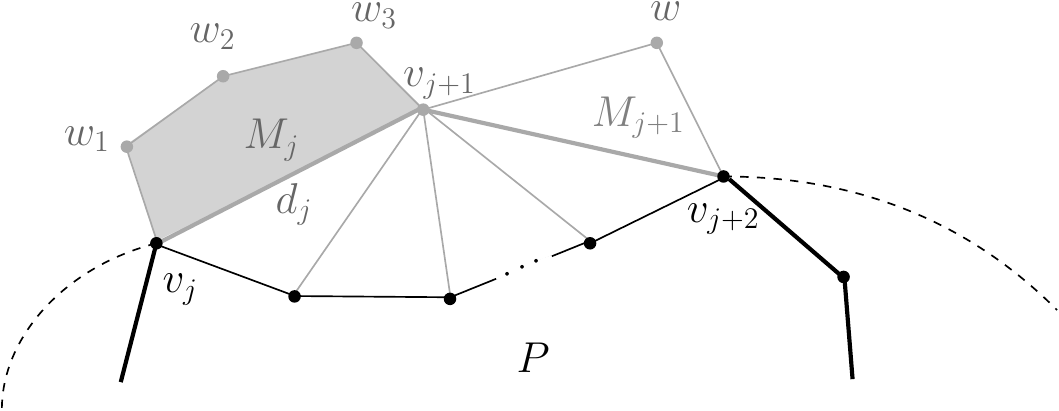}
		\label{fig:case7b}
	}\\
	\subfloat[$D''=\{v_j, v_{j+1}\}$ and $v_j\in D'$.]{
		\includegraphics[scale=0.46,page=34]{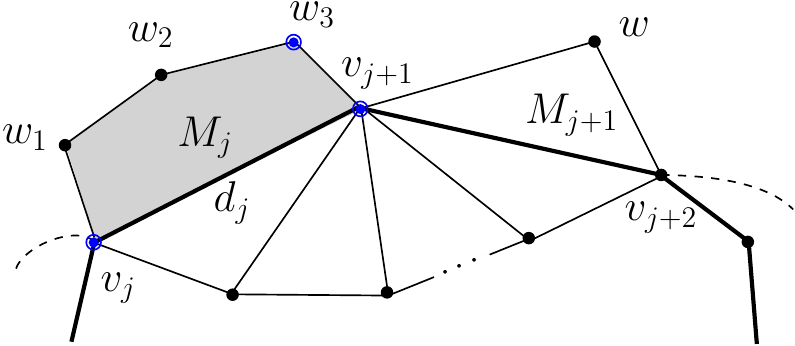}
		\label{fig:case7_1}
   }~~~~
   	\subfloat[$D''=\{v_{j+1}, w_2\}$ and $v_j\notin D''$.]{
		\includegraphics[scale=0.46,page=35]{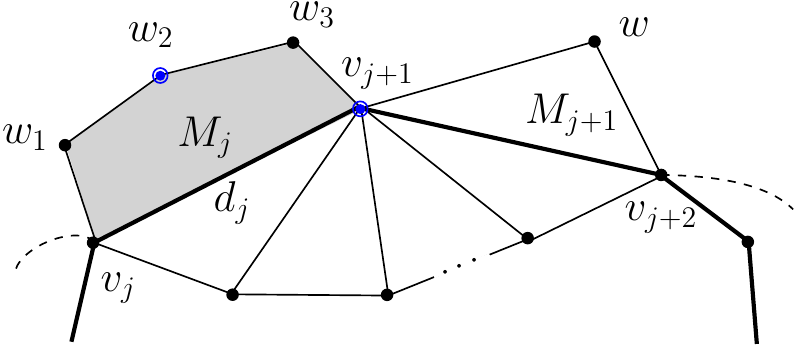}
		\label{fig:case7_2}
   }
    \caption{Illustrating Case 7 of Theorem~\ref{the:boundpaired}.}
    \label{fig:case7}
\end{figure}

Suppose that there exist two consecutive MOP graphs $M_j$ and $M_{j+1}$ of orders 5 and 3, respectively (the same reasoning can be applied if $|M_j| = 3$ and $|M_{j+1}| = 5$). Let $v_j, w_1, w_2, w_3, v_{j+1}$ and $v_{j+1}, w, v_{j+2}$ be the vertices of $M_j$ and $M_{j+1}$, respectively. 
We note that the following analysis can be performed regardless of whether $d_j$ is contractible or not.

Consider the near-triangulation $T_0$ of order $n-4$ obtained by removing $w_1, w_2, w_3$ and $w$ from $T$. Observe that $T_0$ contains the terminal polygon $P$, which is at least a triangle with an interior vertex, and the near-triangulation $M_k$, which is at least a triangle. Thus, $n-4\ge 5$.
%


Since $P$ has no diagonals and at least one interior vertex, by applying Lemma~\ref{lem:RemovingPoints}~$i)$ to $P$, we derive that the vertex $v_{j+1}$ can be removed from $T_0$ to obtain a near-triangulation $T'$ of order $n-5\ge 4$ (see Figure~\ref{fig:case7b}). By the inductive hypothesis, $T'$ has a PD-set $D'$ of size at most $2\left \lfloor \frac{n-5}4 \right \rfloor$
and, by Lemma~\ref{lem:pentagon}, there is a TD-set $D''$ in $M_j$ of size 2 containing $v_{j+1}$.

We construct a PD-set $D$ in $T$ from $D'$ and $D''$ as follows. If $v_j\in D'\cap D''$, then $D=D'\cup\{v_{j+1},w_3\}$ is a PD-set in $T$ of size at most $2\left \lfloor \frac{n-5}4 \right \rfloor+2 \le \cota$, where $w_3$ and $v_{j+1}$ are partners in $D$; see Figure~\ref{fig:case7_1}. If $v_j\notin D'\cap D''$, then $D'\cap D'' = \emptyset$, so $D=D'\cup D''$ is a PD-set in $T$ of size at most $2\left \lfloor \frac{n-5}4 \right \rfloor+2\le \cota$; see Figure~\ref{fig:case7_2}.


\vskip 0.2 cm
\noindent
{\bf Case 8:} $|M_j|=|M_{j+1}|=5$, with $d_j$ and $d_{j+1}$
not contractible  in $T_j$ and $T_{j+1}$, respectively,  for some $j\in \{1,\dots,k-2\}$.
\vskip 0.1cm

Assume that none of the preceding cases holds. In particular, this implies that if $M_k$ is a MOP graph, then $|M_k| = 5$, and if $M_k$ is not a MOP graph, then $|M_k| \ge 6$ as $T$ is irreducible.

\begin{figure}[tb]
	\centering
	\subfloat[The near-triangulation $T'$ obtained from $T$.]{
		\includegraphics[scale=0.46,page=10]{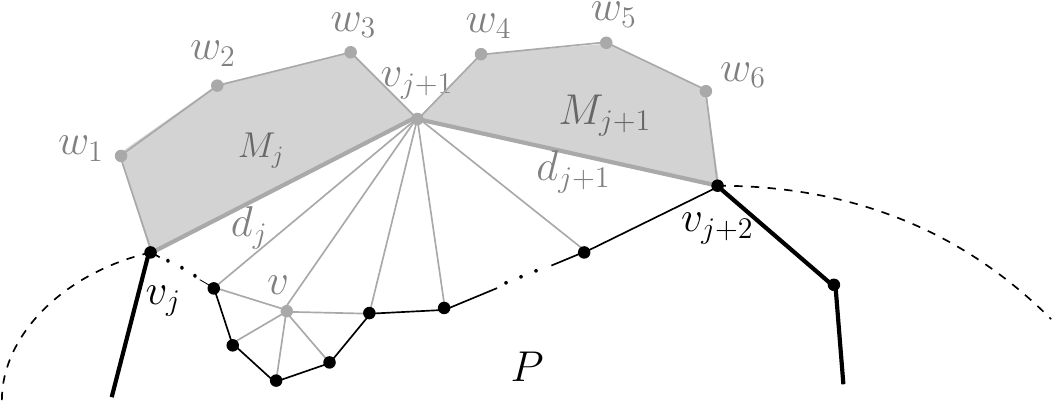}
		\label{fig:case8b}
	}\\
	\subfloat[$D''=\{v_j, v_{j+1}\}$ and $v_j\in D'$.]{
		\includegraphics[scale=0.46,page=36]{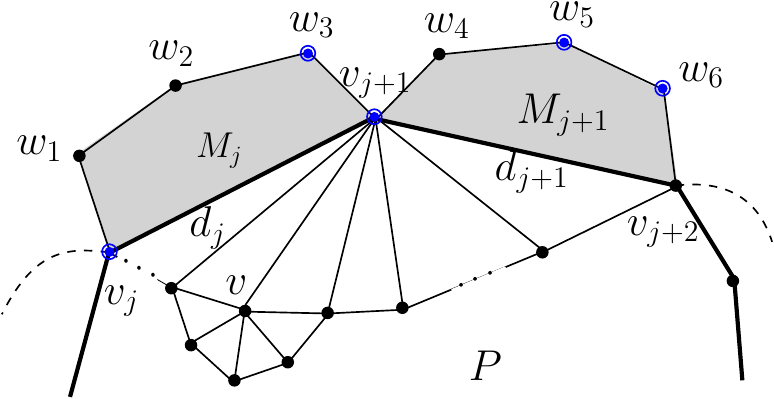}
		\label{fig:case8_1}
   }~~~~~~
   	\subfloat[$D''=\{v_{j+1}, w_2\}$ and $v_j\notin D''$.]{
		\includegraphics[scale=0.46,page=37]{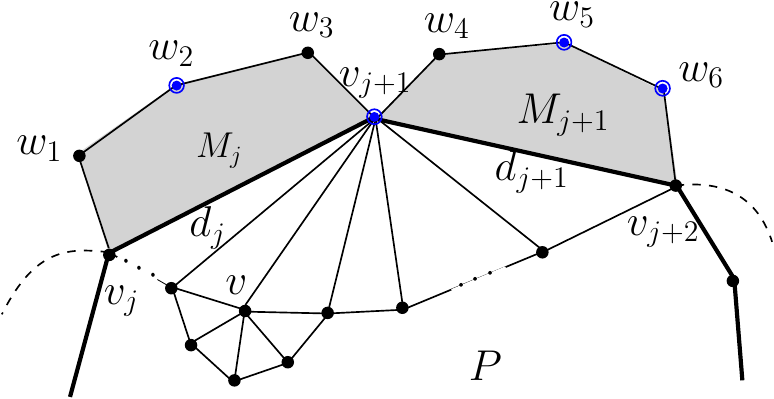}
		\label{fig:case8_2}
   }
    \caption{Illustrating Case 8 of Theorem~\ref{the:boundpaired}.}
    \label{fig:case8}
\end{figure}

Let $M_j$ and $M_{j+1}$ be MOP graphs of order $5$ such that $d_j$ and $d_{j+1}$ are not contractible in $T_j$ and $T_{j+1}$, respectively.
Let $v_j, w_1, w_2, w_3, v_{j+1}$ and $v_{j+1},$ $w_4,$ $w_5,$  $w_6,$ $v_{j+2}$ be the vertices of $M_j$ and $M_{j+1}$, respectively.  
By removing $w_1, w_2, w_3, w_4, w_5$ and $w_6$ from $T$, we obtain a near-triangulation $T_0$ of order $n-6$.

Observe that $n-6 \ge 7$ because $T_0$ contains $P$ and $M_k$, where $P$ is at least a triangle with at least one interior vertex and $M_k$ has order at least 5.
Furthermore, since $d_j$ is not contractible in $T_j$, it is not contractible in either $P$ or $T_0$. Thus, as $P$ has no diagonals, by Lemma~\ref{lem:RemovingPoints} $ii)$ there is an interior vertex $v$ in $P$ (in $T_0$), adjacent to $v_{j+1}$, such that $T'=T_0 - \{v_{j+1},v\}$ is a near-triangulation of order $n-8 \ge 5$ (see Figure~\ref{fig:case8b}).

By the inductive hypothesis, $T'$ has a PD-set $D'$ of size at most $2\left \lfloor \frac{n-8}4 \right \rfloor$ and, by Lemma~\ref{lem:pentagon}, there is a TD-set $D''$ in $M_j$ of size 2 containing $v_{j+1}$.
We construct a PD-set $D$ in $T$ from $D'$ and $D''$ as follows. If $v_j\in D'\cap D''$, then $D=D'\cup\{v_{j+1},w_3\}\cup\{w_5,w_6\}$ is a PD-set in $T$ of size at most $2\left \lfloor \frac{n-8}4 \right \rfloor+4= \cota$, where $w_3$ is paired with $v_{j+1}$, and $w_5$ is paired with $w_6$ in $D$;  see Figure~\ref{fig:case8_1}. If $v_j\notin D'\cap D''$, then $D=D'\cup D''\cup\{w_5,w_6\}$ is a PD-set in $T$ of size at most $2\left \lfloor \frac{n-8}4 \right \rfloor+4= \cota$;  see Figure~\ref{fig:case8_2}.

Therefore, in all cases, we have been able to construct a paired dominating set in $T$ of size at most $\cota$. Consequently, the inductive step is complete, and the theorem follows.
\end{proof}

\section{Semipaired domination in near-triangulations}\label{sec:semipaired}

This section is devoted to showing that $\gamma_{pr2}(T) \le  \left \lfloor \frac{2n}{5} \right \rfloor$ for any near-triangulation $T$ of order $n\ge 5$, with a few exceptions.  Throughout the section, we use $d_T(x,y)$ to denote the distance between two vertices $x$ and $y$ in a near-triangulation $T$. In addition, if $d_T(x,y) \le 2$ we say that two vertices $x$ and $y$ \emph{can be semipaired} in a near-triangulation $T$.

The proof of this upper bound follows the lines of the proof of the preceding theorem and the proof given in~\cite{Claverol21}. Of course, the proof must be adapted to the new scenario of semipaired domination. In relation to the proof given in Theorem~\ref{the:boundpaired}, the main differences are the following: (1) we have to analyze different cases, since in general we can add two vertices to a semi-PD-set for every 5 vertices when applying the inductive hypothesis, (2) there are some near-triangulations for which the upper bound $\left \lfloor \frac{2n}{5} \right \rfloor$ does not hold, hence we have to take into account these exceptions in the inductive process, (3) when obtaining a near-triangulation $T'$ from a near-triangulation $T$ to apply induction on $T'$, semipaired vertices that are at distance 2 in $T'$ can be at distance 3 in $T$, so we have to analyze carefully these situations to ensure that the distance between two semipaired vertices is at most 2 in $T$.

\begin{figure}[t!]
	\centering
	\includegraphics[scale=0.65,page=38]{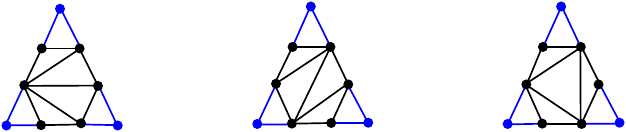}
	\caption{Up to isomorphism, the MOP graphs that belong to $\mathcal{F}$.}
	\label{fig:efe}
\end{figure}

\begin{figure}[t!]
	\centering
	\includegraphics[scale=0.65,page=11]{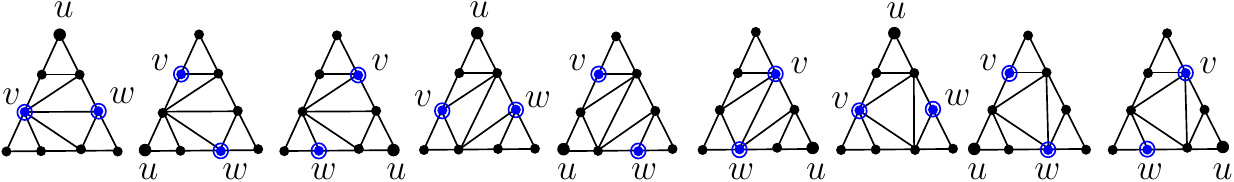}
	\caption{In blue, two vertices $v$ and $w$ at distance at most 2 dominating $V\setminus \{u\}$, for every MOP graph in $\mathcal{F}$.}
	\label{fig:SPDSinF}
\end{figure}

The bound $\left \lfloor \frac{2n}{5} \right \rfloor$ for MOP graphs was already proved in~\cite{Henning19}. Let $\mathcal{F}$ be the family of MOP graphs of order 9 defined as follows~\cite{Henning19} (see Figure~\ref{fig:efe}): Add a vertex adjacent to the endvertices of each of three alternating edges on the boundary of the outer face of a MOP graph of order 6. 
Henning and Kaemawichanurat proved the following theorem.

\begin{theorem}[\cite{Henning19}]\label{the:henning}
	If $T$ is a MOP graph of order $n\ge 5$ not belonging to $\mathcal{F}$, then $\gamma_{pr2}(T) \le \left \lfloor \frac{2n}{5} \right \rfloor$ and the bound is tight.
\end{theorem}

Thus, MOP graphs in $\mathcal{F}$ do not have semi-PD-sets of cardinality 2. As these MOP graphs will appear in the proof of Theorem~\ref{the:boundsemipaired} when applying induction, the following two lemmas show dominating sets or semi-PD-sets in some cases.

\begin{lemma}\label{casiTDSmop9} If $T=(V,E)$ is a MOP graph in $\mathcal{F}$ and $u$ is a vertex of degree 2 in $T$, then $T$ has a set $D=\{v,w\}$ of cardinality 2 dominating every vertex in $V\setminus \{u\}$ and such that $d_T(v,w)\le 2$, $d_T(u,v)= 2$ and $d_T(u,w)= 2$.
\end{lemma}
\begin{proof}
	A set satisfying the given conditions is illustrated in Figure~\ref{fig:SPDSinF} for all MOP graphs in $\mathcal{F}$.
\end{proof}

Observe that in every case of Figure~\ref{fig:SPDSinF}, by adding an edge connecting $u$ to some of $v$ and $w$, the vertices of $D$ in Lemma~\ref{casiTDSmop9}, we obtain a near-triangulation of order 9 that contains a semi-PD-set of size 2. With this observation, the following lemma holds.

\begin{lemma}\label{semiPDset9} Let $T$ be a near-triangulation of order 9 such that when removing a boundary edge $e$, the resulting near-triangulation is a MOP graph in $\mathcal{F}$. Then $T$ contains a semi-PD-set of size 2, so $\gamma_{pr2}(T) = \left \lfloor \frac{2n}{5} \right \rfloor$.
\end{lemma}
\begin{proof}
	If we obtain a MOP graph $H$ in $\mathcal{F}$ after removing $e$, then necessarily one of the endpoints of $e$ is a vertex $u$ of degree two in $H$, and the other is a vertex at distance two from $u$ in $H$. Thus, by the previous observation, the set $D$ of Lemma~\ref{casiTDSmop9} associated with $u$ is a semi-PD-set of size 2 in $T$.
\end{proof}

The following technical lemma will be necessary in the proof of the main theorem to ensure that the elements of a dominating set can be semipaired 
into pairs of vertices at distance at most 2.

\begin{lemma}\label{lem:recombinandoparejas}
	Let $e=uv$ be a contractible edge of a near-triangulation $T$.
	If $\{(x_i,y_i):i=1,\dots ,r\}$ is a set of $r$ pairs of vertices of $T$ different from $u$ and $v$ such that $d_{T/e}(x_i,y_i)= 2$ and $d_T(x_i,y_i)=3$,
    then there is a partition of $\{x_i:1\le i\le r\}\cup \{y_i:1\le i\le r\}$ into $r$ pairs such that, with at most one exception, the vertices of each pair are both in $N_T(u)$ or both in $N_T(v)$.

\end{lemma}
\begin{proof}
	Denote by $w$ the vertex of $T/e$ corresponding to the edge $e$, and by $A$ the set  $\{x_i:1\le i\le r\}\cup \{y_i:1\le i\le r\}$.
	If  $d_{T/e}(x_i,y_i)= 2$ and $d_T(x_i,y_i)=3$,  then  $w$ is the only common neighbor of $x_i$ and $y_i$ in $T/e$.
	Hence,  one of the vertices of $\{x_i,y_i\}$, say $x_i$, is a neighbor of $u$, and then $y_i$ is a neighbor of $v$.
	Thus, $|A\cap N_T(u)|=|A\cap N_T(v)|=r$.
	If $r$ is even, then consider any partition of $A\cap N_T(u)$ into 
    pairs and any partition of $A\cap N_T(v)$ into 
    pairs to obtain the desired result.
	If $r$ is odd, consider a partition of $A$ with 
    pair formed by an element in  $A\cap N_T(u)$ and another in $A\cap N_T(v)$, together with a partition of the remaining vertices of $A\cap N_T(u)$ into 
    pairs and the remaining vertices of $A\cap N_T(v)$ into 
    pairs..
\end{proof}

Now we are ready to prove the main theorem of this section.

\begin{theorem}\label{the:boundsemipaired}
	If $T=(V,E)$ is a near-triangulation of order $n\ge 5$ not belonging to $\mathcal{F}$, then $\gamma _{pr2}(T)\le  \left \lfloor \frac{2n}{5} \right \rfloor$ and the bound is tight.
\end{theorem}

\begin{proof}
As before, we proceed by induction on the order and the number of interior vertices of near-triangulations. Notice that since the bound is tight for MOP graphs (Theorem~\ref{the:henning}), it is also tight for near-triangulations.

Let $T$ be a near-triangulation of order $n\ge 5$, not belonging to $\mathcal{F}$, with $m\ge 0$ interior vertices. If $m=0$, then $T$ is a MOP graph, and the result follows by Theorem~\ref{the:henning}. Furthermore, if $n\in \{5,6\}$ and $m>0$, then $T$ is reducible. Thus, by successively removing boundary edges from $T$, we obtain a MOP graph $T'$ of order $n$,  $n\in \{5,6\}$, such that $\gamma_{pr2}(T)\le \gamma_{pr2} (T')\le \cotas$.
	
Now, let $T$ be a near-triangulation of order $n\ge 7$, with $m>0$ interior vertices.
We will prove that the result holds for $T$, assuming that the bound holds for every near-triangulation $T'$ not in $\mathcal{F}$ of order $n'$ with $m'$ interior vertices such that either $n'<n$, or $n=n'$ and $m'<m$.
	
Suppose that $T$ is reducible. We can remove a suitable boundary edge $e$ to obtain a near-triangulation $T'$ of order $n$ with $m-1$ interior vertices. If $T'$ does not belong to $\mathcal{F}$, by the inductive hypothesis $\gamma _{pr2}(T') \le \left \lfloor \frac {2n}5 \right \rfloor$, so $\gamma _{pr2}(T) \le \gamma _{pr2}(T') \le \cotas$. Otherwise, $T'$ belongs to $\mathcal{F}$, so $\gamma _{pr2}(T) = \cotas$ by Lemma~\ref{semiPDset9}. In any case, the result holds.

Suppose now that $T$ is irreducible. Then, recall that $T$ contains a terminal polygon $P$ with $k\ge 3$ sides, $d_1=v_1v_2$, $d_2=v_2v_3$, $\dots $, $d_k=v_kv_1$, and that for $j=1,\ldots ,k$, we have $T_j = T_{in}(P,d_j)$; $M_j = T_{out}(P,d_j)$, where $M_j$ is a MOP graph except possibly $M_k$; $|T_j| \ge 6$, and $T_j-d_j$ is a near-triangulation.

We distinguish cases depending on the order of the MOP graphs $M_1, M_2, \ldots , M_k$. To apply induction, we only need to remove the vertices of one MOP graph if there is a MOP graph $M_j$ of order 4 or at least 6 (Cases 1-4), and to remove the vertices from two consecutive MOP graphs if all MOP graphs $M_j$ have order 3 or 5 (Cases 5-7).

	
\vskip 0.2 cm
\noindent
{\bf Case 1:} $|M_j|=4$, for some $j\in \{1,\dots, k-1\}$.
\vskip 0.1 cm
	
Suppose that there is a MOP graph $M_j$ of order 4, whose vertices are $v_j, w_1, w_2$ and $v_{j+1}$ in clockwise order. We build a near-triangulation $T'$ of order $n$, with $m-1$ interior vertices, as in Case 1 of Theorem~\ref{the:boundpaired} (see Figure~\ref{fig:case1-size4s}). Without loss of generality, we may assume that $v_j$ has degree 3 in $M_j$. We remove $w_1$, $w_2$ and $d_j$, and we add $x$ connected to $v_j$ and $v$, and $y$ connected to $x$ and $v_j$, where $v$ is the boundary vertex in $T$ preceding $v_j$. In such a case, we claim that $T'$ does not belong to $\mathcal{F}$. Indeed, consider the third vertex $w$ of the triangular face of $P$ that has $d_j$ as a boundary edge. As $w$ is an interior vertex in $T$, its degree is at least 3 in both $T$ and $T'$. Thus, the vertices $v_j$, $w$, and $v_{j+1}$ are three consecutive vertices on the boundary of $T'$ having degree at least 3. Hence, $T'\notin \mathcal{F}$, as claimed.

\begin{figure}[!htb]
	\centering
    \subfloat[The near-triangulation $T'$ obtained from $T$.]{
	\includegraphics[scale=0.5,page=39]{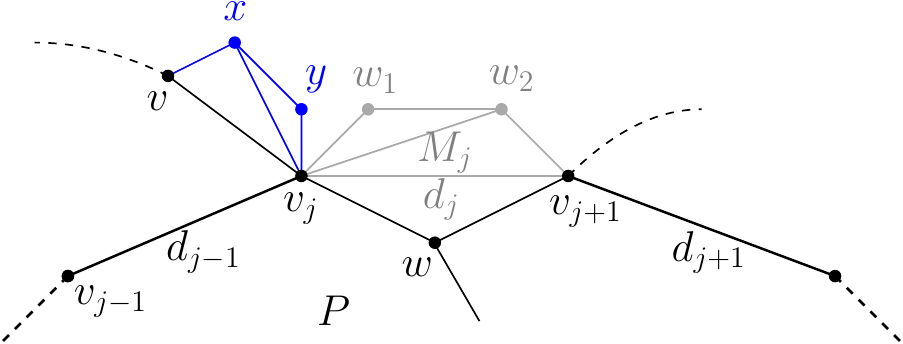}
	\label{fig:case1-size4s}
	}\\
	\subfloat[$v_{j}\in D'$ and $\{x,y\}\notin D'$.]{
		\includegraphics[scale=0.46,page=40]{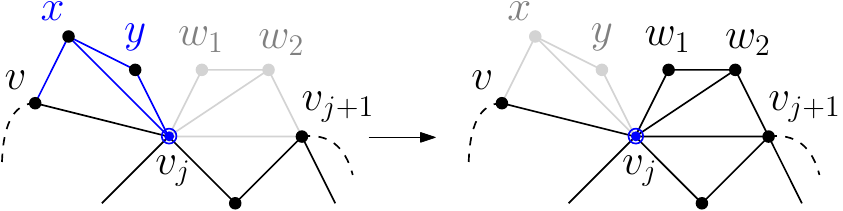}
		\label{fig:m4sa}
	}~~~~
	\subfloat[$v_{j}\in D'$, $x\notin D'$ and $y\in D'$.]{
		\includegraphics[scale=0.46,page=41]{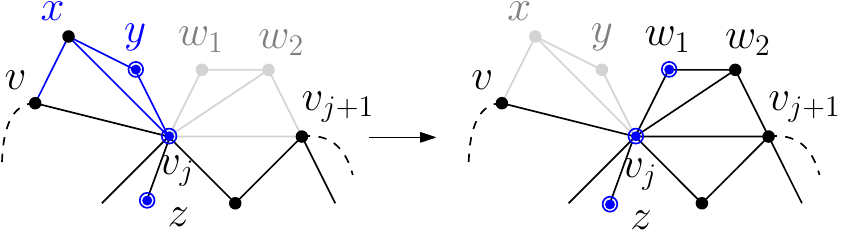}
		\label{fig:m4sb}
   }\\
   	\subfloat[$x$ is semipaired with $z\in N_{T'}{[v_j]}$.]{
		\includegraphics[scale=0.46,page=42]{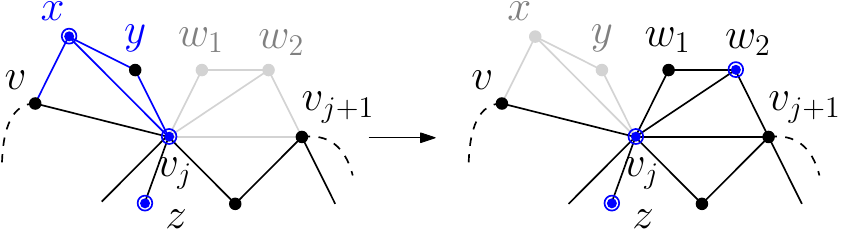}
		\label{fig:m4sc}
	}~~~~
	\subfloat[$x$ is semipaired with $z\in N_{T'}(v)$ and $x'\notin D'$.]{
		\includegraphics[scale=0.46,page=43]{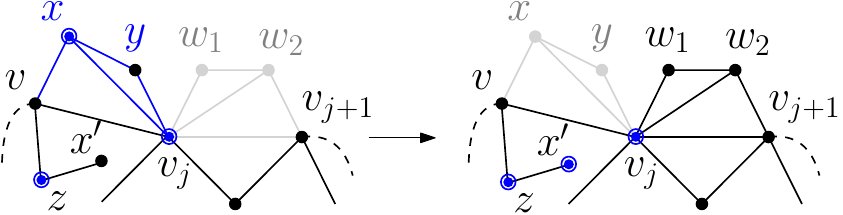}
		\label{fig:m4sd}
   }\\
   	\subfloat[$y$ is semipaired with $z$, $x$ with $z'$, and $d_{T}(z,z')\le 2$.]{
		\includegraphics[scale=0.46,page=44]{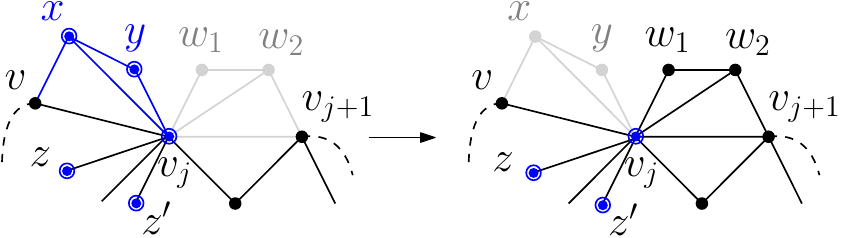}
		\label{fig:m4se}
	}~~~~
	\subfloat[$y$ is semipaired with $z\in N_{T'}(v_j)$, $x$ with $z\in N_{T'}(v)$, and $x'\notin D'$.]{
		\includegraphics[scale=0.46,page=45]{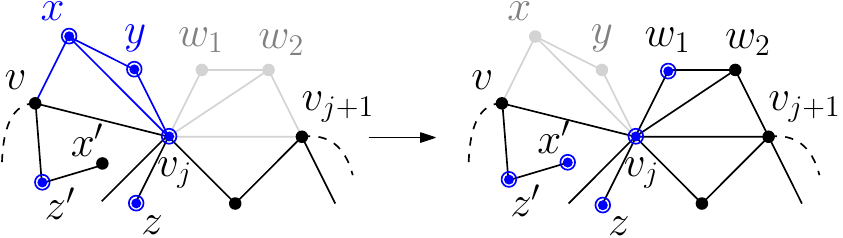}
		\label{fig:m4sf}
   }\\
      	\subfloat[$v_{j}\notin D'$, $x\in D'$ and $y\notin D'$.]{
		\includegraphics[scale=0.46,page=46]{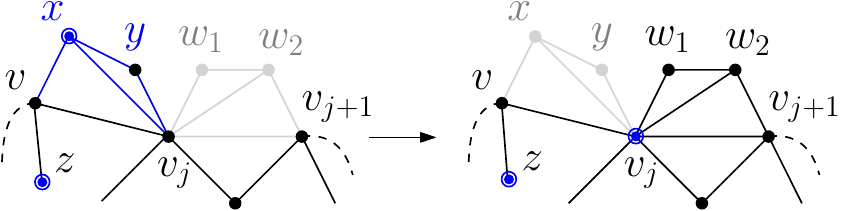}
		\label{fig:m4sg}
	}~~~~
	\subfloat[$v_{j}\notin D'$, $x\notin D'$ and $y\in D'$.]{
		\includegraphics[scale=0.46,page=47]{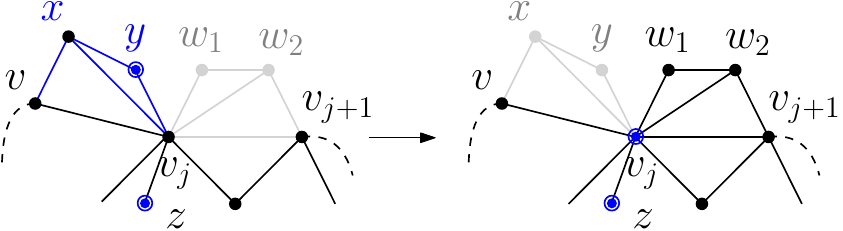}
		\label{fig:m4sh}
   }\\
   	\subfloat[$v_{j}\notin D'$, and $\{x, y\}\in D'$.]{
		\includegraphics[scale=0.46,page=48]{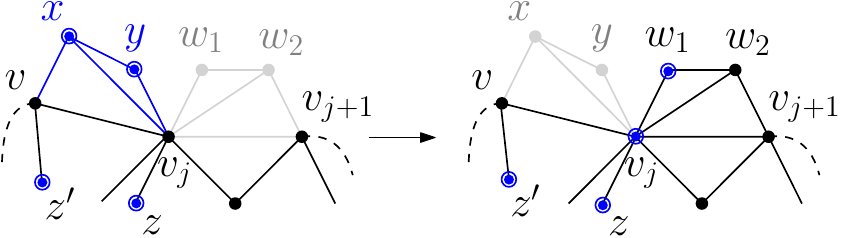}
		\label{fig:m4si}
    }
    \caption{Illustrating Case 1 of Theorem~\ref{the:boundsemipaired}.}
    \label{fig:m4s}
\end{figure}

By the inductive hypothesis, there is a semi-PD-set $D'$ in $T'$ of size at most $\cotas $. From $D'$, we can obtain a semi-PD-set $D$  in $T$ of size at most $\cotas $ as follows. Observe first that $d_{T}(u,v) \le d_{T'}(u,v)$ for any two vertices $u$ and $v$ different from $x$, $y$, $w_1$ and $w_2$, so if two vertices different from $x$ and $y$ are semipaired in $D'$, they can also be semipaired in $T$. We now analyze several cases, depending on whether $v_{j}, x$ and $y$ belong to $D'$ or not.

Suppose first that $v_{j}$ belongs to $D'$. 

\begin{enumerate}
	\item[$\bullet$]  If $\{x,y\}\notin D'$, then $D=D'$ is a semi-PD-set in $T$, since $w_1,w_2\in N_T(v_j)$; see Figure~\ref{fig:m4sa}.
			
	\item[$\bullet$]  If  $x\notin D'$ and $y\in D'$, let $z$ be the vertex in $D'$ semipaired with $y$ in $D'$; see Figure~\ref{fig:m4sb} for an illustration. Thus, $z\in N_{T'}[v_j]$. 
    Since $d_T(z,w_1)\le 2$, the set $D=(D'\setminus \{y\}) \cup \{w_1\}$ is a semi-PD-set in $T$, where $w_1$ is semipaired with $z$ in $D$.
			
	\item[$\bullet$]  If $x\in D'$ and $y\notin D'$, let $z$ be the vertex semipaired with $x$ in $D'$. Thus, $z\in N_{T'}[v]\cup N_{T'}[v_j]$.
			
			\begin{enumerate}
				\item[$-$]  If $z\in N_{T'}[v_j]$, then $D=(D'\setminus\{x\})\cup \{w_2\}$ is a semi-PD-set in $T$, where $w_2$ is semipaired with $z$ in $D$; see Figure~\ref{fig:m4sc}. 
				\item[$-$] If $z\notin N_{T'}[v_j]$,  it suffices to analyze the case $z\in N_{T'}(v)$, since $v\in N_{T'}[v_j]$. In such a case, if $N_{T'}(z)\subseteq D'$, that is, all the neighbors of $z$ belong to $D'$, then $D=D'\setminus\{x,z\}$ is a semi-PD-set in $T$.
				Otherwise, there is a neighbor $x'$ of $z$ not in $D'$, so $D=(D'\setminus\{x\})\cup \{x'\}$ is a semi-PD-set in $T$, with $x'$ semipaired with $z$ in $D$; see Figure~\ref{fig:m4sd}.
			\end{enumerate}
			
		\item[$\bullet$] If $x\in D'$ and $y\in D'$, let $z$ be the vertex semipaired with  $y$ and $z'$ the vertex semipaired with $x$ in $D'$. Thus, $z\in N_{T'}[v_j]$ and $z'\in N_{T'}[v]\cup N _{T'}[v_j]$.

      \begin{enumerate}
\item[$-$] If $d_T(z,z')\le 2$, then $D=D'\setminus \{x,y\}$ is a semi-PD-set in $T$, where  $z$ is semipaired with $z'$ in $D$; see Figure~\ref{fig:m4se}.
			
\item[$-$] If $d_T(z,z')=3$, then $z\in N_{T'}(v_j)\setminus N _{T'}[v]$ and $z'\in N_{T'}(v)\setminus N _{T'}[v_j]$. In this case, if $N_{T'}(z')\subseteq D'$, then $D=(D'\setminus\{x,y,z'\})\cup \{w_1\}$ is a semi-PD-set in $T$ with $w_1$ semipaired with $z$. If there is a neighbor $x'$ of $z'$ not in $D'$, then  $D=(D'\setminus\{x,y\})\cup \{w_1,x'\}$ is a semi-PD-set in $T$ with $w_1$ semipaired with $z$, and $x'$ semipaired with $z'$; see Figure~\ref{fig:m4sf}.
\end{enumerate}

\end{enumerate}

Now suppose that $v_{j}\notin D'$. 		

		\begin{enumerate}
			\item[$\bullet$]
			If $y\notin D'$, then $x$ must belong to $D'$ to dominate $y$; see Figure~\ref{fig:m4sg}. Let  $z\in (N_{T'}[v]\setminus \{v_{j}\})\cup N_{T'}(v_j )$ be the vertex semipaired with $x$ in $D'$. In this case, $D=(D'\setminus \{x\}) \cup \{v_{j}\}$ is a semi-PD-set in $T$, where $z$ and $v_j$ are semipaired in $D$.
			
			\item[$\bullet$]
			If $y\in D'$ and $x\notin D'$, let  $z\in N_{T'}(v_j)$ be the vertex semipaired with $y$ in $D'$; see Figure~\ref{fig:m4sh}. Thus, $D=(D'\setminus \{y\}) \cup \{v_j\}$ is a semi-PD-set in $T$, semipairing $z$ and $v_j$ in $D$.
			
			\item[$\bullet$] If $y\in D'$ and $x\in D'$,  let $z$ be the vertex semipaired with  $y$ and $z'$ be the vertex semipaired with $x$ in $D'$; see Figure~\ref{fig:m4si}. Thus, $z\in N_{T'}(v_j)$ and $z'\in (N_{T'}[v]\setminus \{v_j\})\cup N _{T'}(v_j)$.			
			Hence, $D=(D'\setminus \{x,y\}) \cup \{v_{j}, w_1\}$ is a semi-PD-set in $T$, where $v_j$ is semipaired with $z'$ and $w_1$ with $z$.
		\end{enumerate}

\vskip 0.2 cm
\noindent
{\bf Case 2:} $|M_j|=6$, for some $j\in \{1,\dots, k-1\}$.
\vskip 0.1 cm
	
Suppose that there is a MOP graph $M_j$ of order 6. Let $v_j, w_1, w_2, w_3,w_4$ and $v_{j+1}$ be the vertices of $M_j$, in clockwise order. By Lemma~\ref{lem:hexagon}, there is a TD-set $D''=\{u,u'\}$ of $M_j$, which contains exactly one of the vertices $v_j$ and $v_{j+1}$, say $v_j=u$, so $u'\neq v_{j+1}$.

Since $v_j$ is adjacent to at least one interior vertex in $P$, by Lemma~\ref{lem:contractibleedge} $i)$, there is an interior vertex $v$ in $P$ adjacent to $v_{j}$, such that the edge $e=vv_j$ is contractible in $P$. Consequently, $e$ is also contractible in $T_j=T-\{w_1,w_2,w_3,w_4\}$. Let $T'=T_j/e$ and let $w$ be the vertex of $T'$ obtained when contracting $e$. Hence, $T'$ is a near-triangulation satisfying $|T'|=|T_j|-1 = n-5\ge 5$.

\begin{figure}[!htb]
	\centering
	\includegraphics[scale=0.40,page=12]{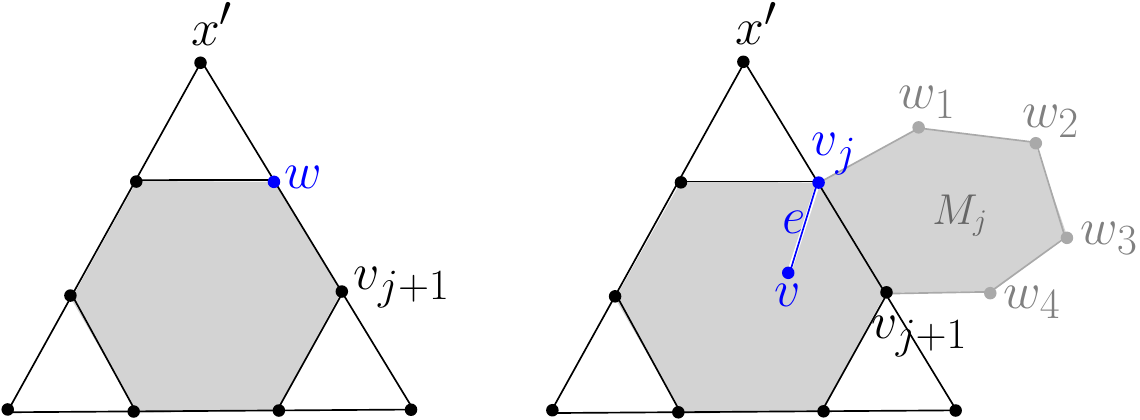}
	\caption{The near-triangulations $T'\in \mathcal{F}$ (left) and $T$ (right). The shaded region is triangulated.} \label{fig:case2F}
\end{figure}

To obtain a semi-PD-set $D'$ in $T'$, we distinguish whether $T'$ belongs to $\mathcal{F}$ or not. If $T'\notin\mathcal{F}$, by the inductive hypothesis, there is a semi-PD-set $D'$ in $T'$ of size at most $\cotasn{n-5}$.

If $T'\in \mathcal{F}$, then $T$ is a near-triangulation of order $14$, with $v$ being the only interior vertex in $P$. Moreover, $w$ is a vertex of degree at least $3$ in $T'$ because $v_j$ is adjacent to at least 3 vertices different from $v$ in $T_j$. Hence, $w$ is adjacent to a vertex $x'$ of degree 2 in $T'$ (see Figure~\ref{fig:case2F}), which must also be adjacent to $v_j$ in $T$ because $v$ is an interior vertex of $P$ in $T$.
By Lemma~\ref{casiTDSmop9}, there exists a set $D'$ consisting of two vertices at distance at most 2 that dominates all vertices of $T'$ except $x'$, such that the distance in $T'$ from $x'$ to either vertex in $D'$ is 2. Note that $D'$ satisfies all the conditions required to be a semi-PD-set in $T'$, except that $D'$ does not dominate $x'$. Moreover, $w\notin D'$  since is adjacent to $x'$ (see Figure~\ref{fig:SPDSinF}), and $|D'| = 2 = \left \lfloor  \frac{2 (14-5)}{5} \right \rfloor = \left \lfloor  \frac{2(n-5)}{5} \right \rfloor $.

\begin{figure}[!htb]
	\centering
	\subfloat[$D''=\{v_j, u'\}$, $x$ is semipaired with $u'$, and $y$ with $v_j$.]{
		\includegraphics[scale=0.48,page=49]{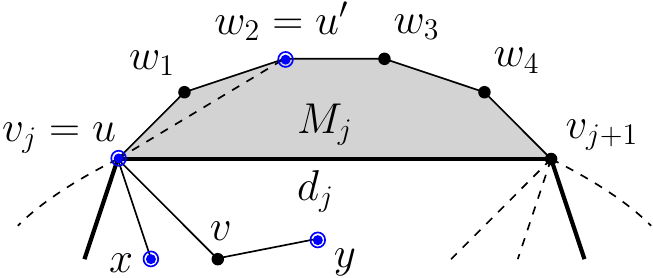}
		\label{fig:ms6a}
	}~~~~~~~~~~~~~~~~
	\subfloat[$D''=\{v_j, u'\}$, $v$ is semipaired with $z$, $x$ with $u'$, and $y$ with $v_j$.]{
		\includegraphics[scale=0.48,page=50]{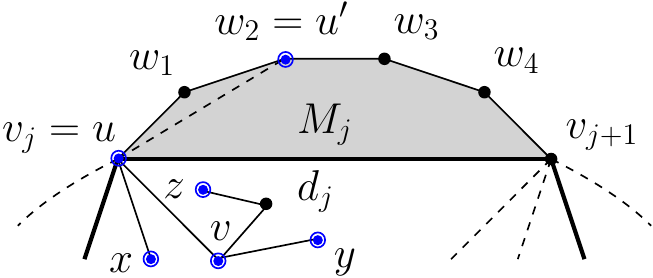}
		\label{fig:ms6b}
    }\\
    \subfloat[$D''=\{v_j, u'\}$, $v$ is semipaired with $y$, $x$ with $u'$, and $z$ with $v_j$.]{
		\includegraphics[scale=0.48,page=51]{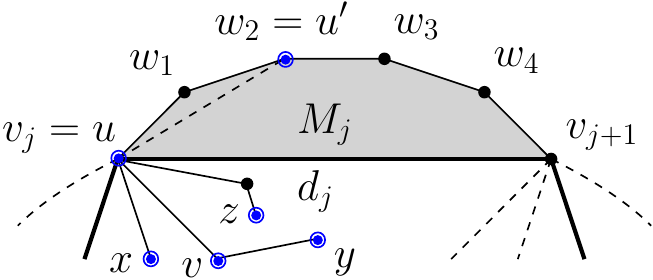}
		\label{fig:ms6c}
    }
    \caption{Illustrating Case 2 of Theorem~\ref{the:boundsemipaired}.}
    \label{fig:ms6}
\end{figure}

Regardless of whether $T'$ belongs to $\mathcal{F}$ or not, from $D'$ and $D''=\{v_j,u'\}$ we can build a semi-PD-set in $T$ as follows. Recall that $D''$ is a TD-set, so $v_j$ and $u'$ are adjacent. We distinguish whether $w$ belongs to $D'$ or not.

\begin{enumerate}
\item[$\bullet$] If $w\notin D'$, we claim that $D=D'\cup D''$ is a semi-PD-set in $T$ of size at most $\cotas$. Clearly, $|D| \le \left \lfloor  \frac{2\ (n-5)}{5} \right \rfloor + 2 \le \cotas$. On the other hand, if two vertices $x$ and $y$ are semipaired in $D'$, they can also be semipaired in $T$ when $d_{T'}(x,y) = d_T(x,y)$, because in such a case $d_{T'}(x,y) \le d_T(x,y)\le d_{T'}(x,y)+1$.
    Consider now the set of $r\ge 0$ pairs $(x_i,y_i)$ in $D'$ such that $x_i$ and $y_i$ are semipaired in $D'$, and $d_{T'}(x_i,y_i) = 2$ and $d_{T}(x_i,y_i) = 3$. If $r>0$, by Lemma~\ref{lem:recombinandoparejas}, we can define a new set of $r$ pairs of vertices such that the two vertices of each pair belong to either $N_T(v_j)$ or $N_T(v)$, so they can be semipaired in $T$, except for at most one pair $(x,y)$. These two vertices $x$ and $y$ are semipaired in $D'$, but they cannot be semipaired in $T$ because $d_T(x,y)=3$, so one of them, say $x$, belongs to $N_T(v_j)$ and $y$ to $N_T(v)$; see Figure~\ref{fig:ms6a}. If such a pair $(x,y)$ does not exist or if $r=0$, $v_j$ is semipaired with $u'$ in $D$, and if such a pair exists, $x$ is semipaired with $u'$ in $D$, and $y$ with $v_j$; see Figure~\ref{fig:ms6a}. In any case, $D$ is a semi-PD-set in $T$ as claimed.
\item[$\bullet$] If $w \in D'$, let $z$ be its semipaired vertex in $D'$. Recall that this case cannot occur if $T'\in \mathcal{F}$ because otherwise $w\notin D'$. We claim that $D=(D'\setminus \{ w\} )\cup D'' \cup \{v \}$ is a semi-PD-set in $T$ of size at most $\cotas$. Observe that $D$ is 
a dominating set in $T$, because $v$ and $v_j$ dominate any vertex that was dominated by $w$ in $T'$. 
Moreover, excluding the pair $(w,z)$, by Lemma~\ref{lem:recombinandoparejas} we may assume that there is at most one pair $(x,y)$ such that $x$ and $y$ are semipaired in $D'$, $x\in N_T(v_j)$, $y\in N_T(v)$ and they cannot be semipaired in $T$. As $z$ is semipaired with $w$ in $D'$, then $d_{T'}(z,w)\le 2$, so $d_T(z,v)\le 2$ or $d_T(z,v_j)\le 2$. Suppose first that $d_T(z,v)\le 2$; see Figure~\ref{fig:ms6b}. If the pair $(x,y)$ exists, then $z$ is semipaired with $v$ in $D$, $y$ with $v_j$ and $x$ with $u'$, and if such a pair does not exist, then $z$ is semipaired with $v$ in $D$, and $v_j$ with $u'$. Suppose now that $d_T(z,v_j)\le 2$; see Figure~\ref{fig:ms6c}. If the pair $(x,y)$ exists, then $z$ is semipaired with $v_j$ in $D$, $y$ with $v$ and $x$ with $u'$, and if such a pair does not exist, then $z$ is semipaired with $v_j$ in $D$, and $v$ with $u'$. In any case, $D$ is a semi-PD-set in $T$, as claimed.
\end{enumerate}



\vskip 2mm
\noindent
{\bf Case 3:} $|M_j|=r\in \{7,8,9\}$, for some $j\in \{1,\dots, k-1\}$.
\vskip 1mm
	
Suppose that there is a MOP graph $M_j$ of order $r$, $r\in \{7,8,9\}$, and consider $T_j$, the near-triangulation obtained from $T$ by removing the vertices of $M_j$, except for the endpoints of $d_j$. Since $T_j$ is a near-triangulation of order $n-r+2$, where $n-r+2\ge 6$, $T_j$ has a semi-PD-set $D'$ of size at most $\cotasn{n-r+2}$, by the inductive hypothesis. Notice that $d_{T_j}(x,y)=d_T(x,y)$, for any two vertices $x$ and $y$ of $T_j$.

Assume first that $M_j\notin \mathcal{F}$. By Theorem~\ref{the:henning}, $M_j$ has a semi-PD-set $D''=\{u,v\}$ of size $2$. We build a semi-PD-set $D$ in $T$ as follows.
If $D'\cap D''=\emptyset$, let  $D=D'\cup D''$.
If $|D'\cap D''|=2$, let  $D=D'$.
If $|D'\cap D''|=1$, suppose $v\in D''\setminus D'$, so $v\notin \{v_j, v_{j+1}\}$. Let $D=D'\cup \{v,z\}$ in this case, where $z$ is a vertex in $M_j$ adjacent to $v$, and different from $v_j$ and $v_{j+1}$, which exists because $r\ge 7$. Clearly, in all cases, $D$ is a semi-PD-set in  $T$, where $v$ is semipaired with $z$ in the last case, satisfying $|D|\le |D'|+2 \le  \cotasn{n-r+2} +2\le \cotas$, when $r\in \{7,8,9\}$.
		
\begin{figure}[!htb]
	\centering
	\includegraphics[scale=0.40,page=13]{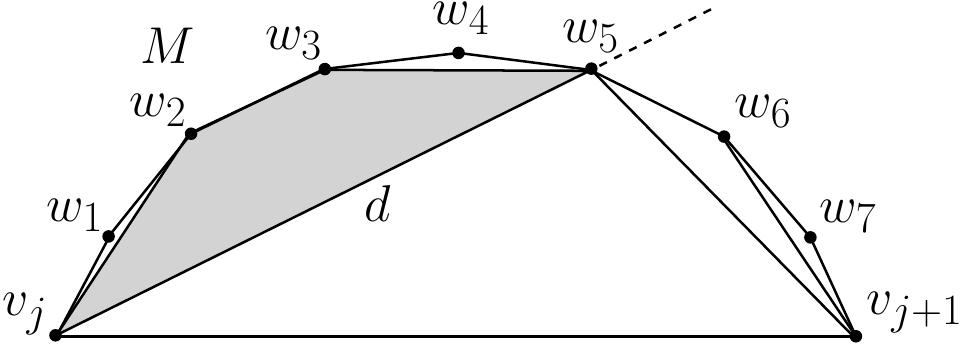}~~~
    \includegraphics[scale=0.40,page=14]{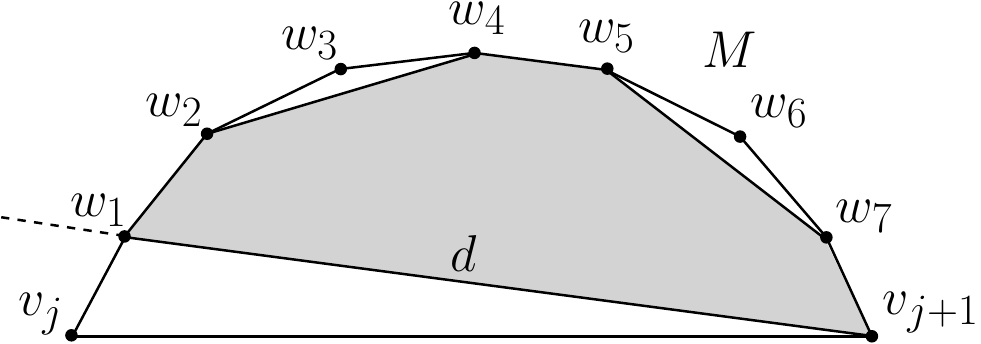}
	\caption{Two examples of the MOP graph $M_j\in \mathcal{F}$. The regions in gray are triangulated. }\label{fig:case9F}
\end{figure}

Assume now that $M_j\in \mathcal{F}$ (see Figure~\ref{fig:efe} for the MOP graphs in $\mathcal{F}$). In such a case, one of the diagonals $d$ of $M_j$ incident to $v_j$ or $v_{j+1}$ splits $M_j$ into two MOP graphs such that the MOP graph $M$ not containing $v_jv_{j+1}$ has order 6, 7 or 8 (see Figure~\ref{fig:case9F} for some examples). If $|M|=8$ or $|M|=7$ (see the right example in Figure~\ref{fig:case9F}), then we proceed as in the two preceding paragraphs, applying the reasoning to $M$ instead of $M_j$.  If $|M|=6$ (see the left example in Figure~\ref{fig:case9F}), consider the near-triangulation $T'$ obtained from $T$ by removing the vertices in $M$ except for the endpoints of $d$. Then, observe that since all neighbors of one of the endpoints of $d$ ($w_5$ in the left example in Figure~\ref{fig:case9F}) are on the boundary, the boundary edge $d$ of $T'$ is contractible in $T'$, by Lemma~\ref{lem:contractibleNew} $ii)$. Thus, we can argue as in Case 2.


\vskip 0.2 cm	
\noindent
{\bf Case 4:} $|M_j|\ge 10$, for some $j\in \{1,\dots, k-1\}$.
\vskip 0.1 cm
	
	Suppose that there is a MOP graph $M_j$ of order $r\ge 10$.  Taking $l=5$ in Lemma~\ref{lem:diagonales}, there exists a diagonal $d$ of $M_j$ dividing $M_j$ into two MOP graphs, such that the MOP graph $M$ not containing $d_j$ has order 6, 7, 8, or 9. Thus, we can proceed as in Cases 2 and 3, where $M$ plays the role of $M_j$.
Recall that if $|M|=6$, then by Lemma~\ref{lem:contractibleNew} $ii)$, $d$ is contractible in the near-triangulation obtained from $T$ by removing the vertices in $M$ except for the endpoints of $d$.

\vskip 2mm
\noindent
{\bf Case 5:} $|M_j|=3$ and $|M_{j+1}|=3$, for some $j\in \{1,\dots,k-2\}$.
\vskip 1mm
	
Suppose that there exist two consecutive MOP graphs $M_j$ and $M_{j+1}$ of order 3. We 
carry out the same construction as in Case 6 of Section~\ref{sec:paired} (see Figure~\ref{fig:case6bs}).
If $v_j, w_1, v_{j+1}$ and $v_{j+1}, w_2, v_{j+2}$ are the vertices of $M_j$ and $M_{j+1}$, respectively, we build a near-triangulation $T'$ of order $n$ with $m-1$ interior vertices, by removing $w_2$, adding a new vertex $w'$ adjacent to $w_1$ and $v_{j+1}$, and removing the edge $d_{j+1}$.

Arguing as in Case 1, since $T'$ contains three consecutive boundary vertices of degree at least 3, namely $v_{j+1}$, $v_{j+2}$ and the boundary vertex $w$ between them, we conclude that $T'\notin \mathcal{F}$. Hence, by the inductive hypothesis, $T'$ has a semi-PD-set $D'$ of cardinality at most $\cotas$.
From $D'$, we obtain a semi-PD-set $D$ in $T$ of size at most $\cotas$, distinguishing two cases. Note that $d_{T}(u,v) \le d_{T'}(u,v)$, for any two vertices $u$ and $v$ different from $w'$ and $w_2$.

\begin{figure}[!htb]
	\centering
	\subfloat[The near-triangulation $T'$ obtained from $T$.]{
		\includegraphics[scale=0.46,page=52]{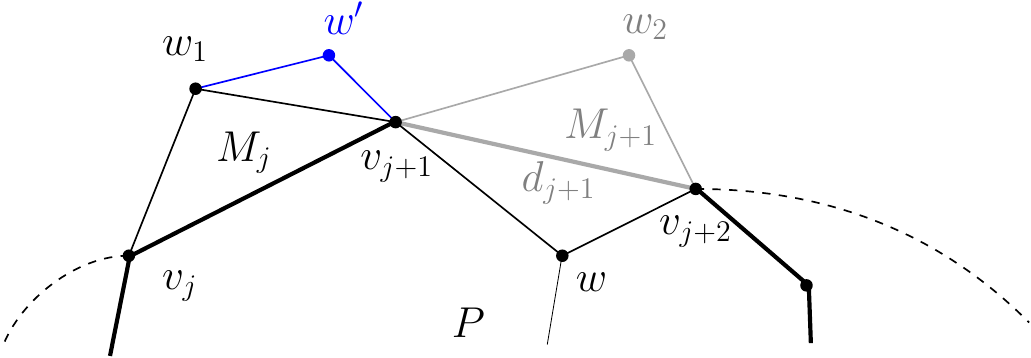}
		\label{fig:case6bs}
	}\\
	\subfloat[$v_{j+1}\in D'$ and $w'\notin D'$.]{
		\includegraphics[scale=0.46,page=53]{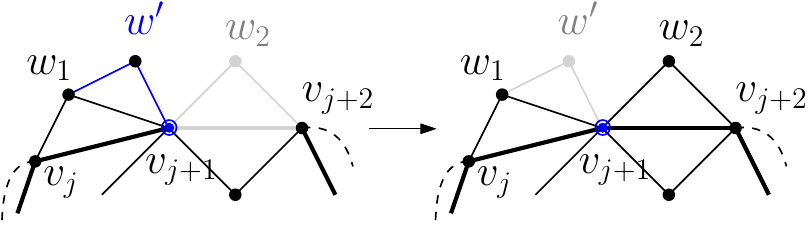}
		\label{fig:cases6_1}
   }~~~~
   	\subfloat[$\{v_{j+1}, w'\}\in D'$ and $w'$ is semipaired with $z$.]{
		\includegraphics[scale=0.46,page=54]{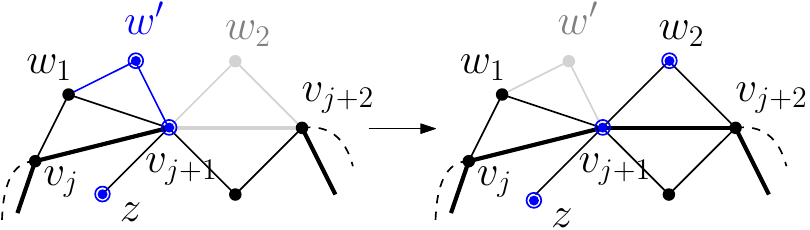}
		\label{fig:cases6_2}
        	}\\
	\subfloat[$\{v_{j+1}, w'\}\notin D'$, and $w_1$ is semipaired with $z$.]{
		\includegraphics[scale=0.46,page=55]{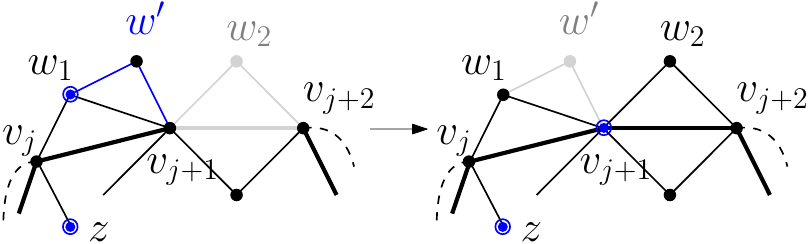}
		\label{fig:cases6_3}
   }~~~~
   	\subfloat[$v_{j+1}\notin D'$, $w'\in D$ and $w'$ is semipaired with $z$.]{
		\includegraphics[scale=0.46,page=56]{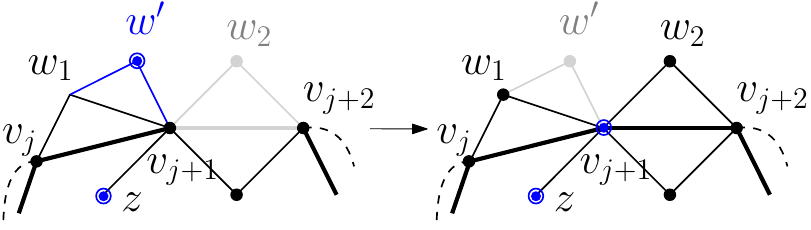}
		\label{fig:cases6_4}
   }
    \caption{Illustrating Case 5 of Theorem~\ref{the:boundsemipaired}.}
    \label{fig:case6s}
\end{figure}

Suppose first that $v_{j+1}\in D'$.		
		\begin{enumerate}
			\item[$\bullet$]
			If $w'\notin D'$, then $D=D'$ is a semi-PD-set in $T$; see Figure~\ref{fig:cases6_1}.
			
			\item[$\bullet$]
			If $w'\in D'$, let $z$ be the semipartner of $w'$ in $D'$. Necessarily, $z\in N_T[v_{j+1}]$; see Figure~\ref{fig:cases6_2}. Then, $D=(D'\setminus \{w'\})\cup \{w_2\}$ is a semi-PD-set in $T$, where $w_2$
			is semipaired with $z$ in $D$.
		\end{enumerate}		

Suppose now that $v_{j+1}\notin D'$. 
				\begin{enumerate}			
			\item[$\bullet$]
			If $w'\notin D'$, then  $w_1$  must belong to $D'$ to dominate $w'$. Let $z$ be the semipartner of $w_1$ in $D'$, that necessarily belongs to $(N(v_{j+1}) \cup N[v_j])\setminus \{v_{j+1}\}$; see Figure~\ref{fig:cases6_3}. Note that since $d_{T'}(w_1,z)\le 2$, then $d_T(v_{j+1},z) \le 2$. Thus,  $D=(D'\setminus \{w_1\})\cup \{v_{j+1}\}$ is a semi-PD-set in $T$, where $v_{j+1}$ is semipaired with $z$ in $D$.
			\item[$\bullet$]
			If $w'\in D'$,  let $z$ be its semipartner in $D'$; see Figure~\ref{fig:cases6_4}. Thus, $z\in N_T(v_{j+1})$, so $D=(D'\setminus \{w'\})\cup \{v_{j+1}\}$ is a semi-PD-set in $T$,
			where $v_{j+1}$ and $z$ are semipaired.
		\end{enumerate}

\vskip 2mm	
\noindent
{\bf Case 6:} $|M_j|=5$ and $|M_{j+1}|=3$, for some $j\in \{1,\dots,k-2\}$.
\vskip 1mm
			
	Suppose that there exist two consecutive MOP graphs $M_j$ and $M_{j+1}$ of orders 5 and 3, respectively. The same reasoning applies if $|M_j|=3$ and $|M_{j+1}|=5$. Let $v_j, w_1, w_2, w_3, v_{j+1}$ and $v_{j+1}, w, v_{j+2}$ be the vertices of $M_j$ and $M_{j+1}$, respectively.
	As in Case 7 of Section~\ref{sec:paired}, consider the near triangulation $T'=T-\{w_1, w_2, w_3,w,v_{j+1}\}$ of order $n-5\ge 4$ (see Figure~\ref{fig:case7bs}). Observe again that $d_{T}(u,v) \le d_{T'}(u,v)$ for any pair of vertices $u$ and $v$ of $T'$.

	If $n-5=4$, then $T$ has order $9$. Consequently, the polygon $P$ must be a triangle with an interior vertex, and there are two consecutive MOP graphs of order $3$ in $T$. Thus, the preceding case applies. Hence, we may assume that $n-5\ge 5$.

\begin{figure}[!httb]
	\centering
	\subfloat[The near-triangulation $T'$ obtained from $T$.]{
		\includegraphics[scale=0.46,page=9]{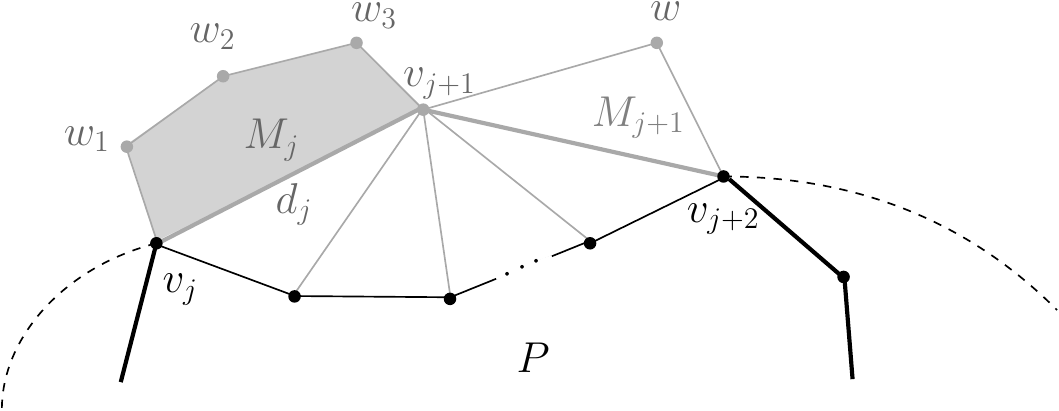}
		\label{fig:case7bs}
	}\\
	\subfloat[$T'\notin \mathcal{F}$; $w_2$ is semipaired with $v_{j+1}$.]{
		\includegraphics[scale=0.46,page=57]{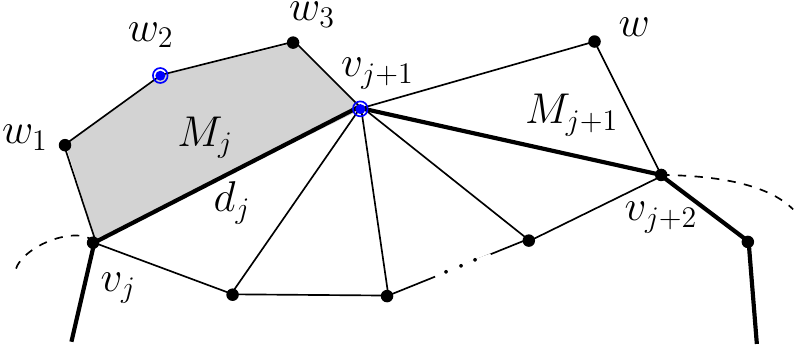}
		\label{fig:cases7_1}
   }~~~~
   	\subfloat[$T'\in \mathcal{F}$; $u$ is semipaired with $v$, and $w_2$ with $v_{j+1}$.]{
		\includegraphics[scale=0.46,page=58]{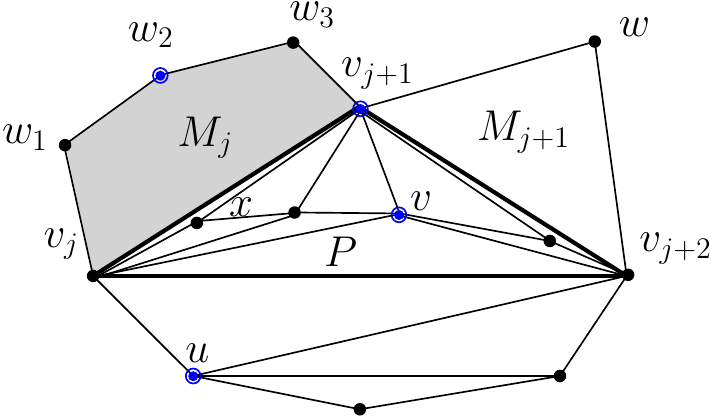}
		\label{fig:cases7_2}
   }
    \caption{Illustrating Case 6 of Theorem~\ref{the:boundsemipaired}.}
    \label{fig:case7s}
\end{figure}

	Suppose first that $T'\notin \mathcal{F}$; see Figure~\ref{fig:cases7_1}. 
    $T'$ has a semi-PD-set $D'$ of size at most $\left \lfloor \frac{2(n-5)}5 \right \rfloor$ by the inductive hypothesis.
	Then, $D=D'\cup \{ v_{j+1}, w_2\}$ is a semi-PD-set in $T$, where $v_{j+1}$ and $w_2$ are semipaired in $D$, and $|D|\le \left \lfloor \frac{2(n-5)}5 \right \rfloor+2\le \cotas$. Note that $d(v_{j+1},w_2)\le 2$ and the set $\{v_{j+1},w_2\}$ dominates all vertices in $M_j$ and $M_{j+1}$.
	
	Suppose now that $T'\in \mathcal{F}$; see Figure~\ref{fig:cases7_2} for an illustration.  Thus, $T$ is a near-triangulation of order $14$. Observe that $v_{j+1}$ is adjacent in $T$ to a vertex $x$ of degree 2 in $T'$, because $v_{j+1}$ has at least 3 neighbors ($v_j$, $v_{j+2}$ and an interior vertex in $P$) in $T-\{w_1, w_2, w_3,w\}$.
	By Lemma~\ref{casiTDSmop9}, there exists a set $D'$ dominating all vertices of $T'$ except possibly $x$, consisting of two vertices at distance at most 2 ($u$ and $v$ in Figure~\ref{fig:cases7_2}).
	Then, $D=D'\cup \{v_{j+1},w_2\}$ is a semi-PD-set in $T$, with $v_{j+1}$ and $w_2$ semipaired in $D$, such that $|D|=4\le \left \lfloor  \frac{2\cdot 14}{5} \right \rfloor = \left \lfloor  \frac{2 n}{5} \right \rfloor $.

\vskip 2mm
\noindent
{\bf Case 7:} $|M_j|=|M_{j+1}|=5$, for some $j\in \{1,\dots,k-2\}$.
\vskip 1mm
	
Assume that none of the preceding cases holds, so all MOP graphs around $P$ have order 5. Let $M_j$ and $M_{j+1}$ be two such MOP graphs. Note that $M_j$ and $M_{j+1}$ are fans, that is, there is a vertex of the boundary, the center of the fan, adjacent to the remaining 4 vertices.
Let $v_j, w_1, w_2, w_3, v_{j+1}$ and $v_{j+1}, w_4, w_5, w_6, v_{j+2}$ be the vertices of $M_j$ and $M_{j+1}$, respectively, in clockwise order.  
By removing $w_1, w_2, w_3, w_4, w_5$ and $w_6$ from $T$, we obtain a near-triangulation $T'$ of order $n-6$; see Figure~\ref{fig:cases8b}.  Observe that since $v_j, v_{j+1}$ and $v_{j+2}$ have degree at least 3 in $T'$, it is a near-triangulation not in $\mathcal{F}$ of order $n-6 \ge 7$. Moreover, $d_{T'}(u,v) = d_T(u,v)$ for any two vertices $u$ and $v$ of $T'$.

\begin{figure}[!htb]
	\centering
	\subfloat[The near-triangulation $T'$ obtained from $T$.]{
		\includegraphics[scale=0.46,page=59]{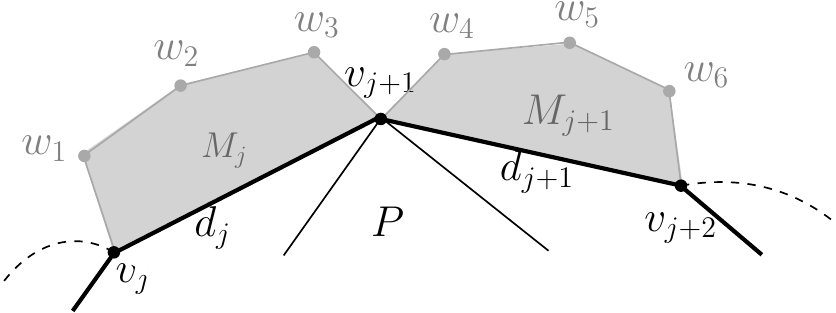}
		\label{fig:cases8b}
	}\\
	\subfloat[$z=z'=v_{j+1}$ and $v_{j+1}\notin D'$.]{
		\includegraphics[scale=0.46,page=60]{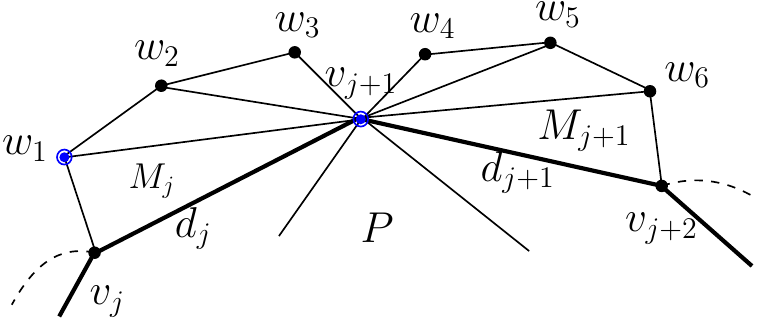}
		\label{fig:cases8_1}
   }~~~~~~
   	\subfloat[$z\neq z'$ and $z'\in D'$.]{
		\includegraphics[scale=0.46,page=61]{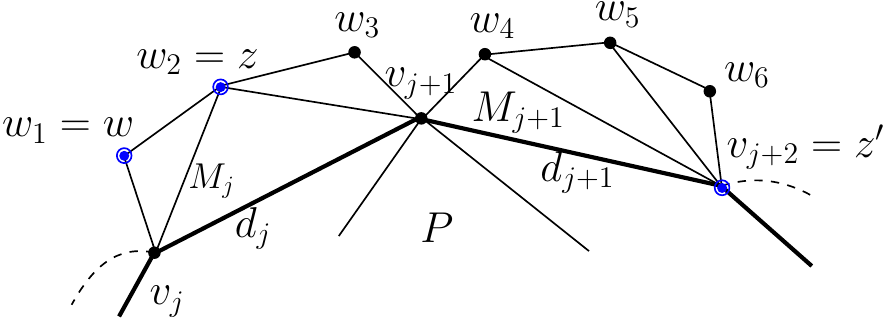}
		\label{fig:cases8_2}
   }
    \caption{Illustrating Case 7 of Theorem~\ref{the:boundsemipaired}.}
    \label{fig:cases8}
\end{figure}

By the inductive hypothesis, $T'$ has a semi-PD-set $D'$ of size at most $\left \lfloor \frac{2(n-6)}5 \right \rfloor$. From $D'$, we construct a semi-PD-set $D$ in $T$ of size $|D|\le |D'|+2\le \cotasn{n-6}+2\le \cotas$ as follows. Let $z$ and $z'$ be the centers of the fans $M_j$ and $M_{j+1}$, respectively. Notice that $d_T(z,z')\le 2$.

\begin{enumerate}
	
\item[$\bullet$]	If $z=z'$, then $z=z'=v_{j+1}$. If $v_{j+1}\in D'$, let $D=D'$, and if $v_{j+1}\notin D'$, let $D=D'\cup \{v_{j+1},w_1\}$, semipairing $v_{j+1}$ and $w_1$ (see Figure~\ref{fig:cases8_1}). In both cases, $D$ is a semi-PD-set in $T$.
	
\item[$\bullet$]	Suppose that $z\not=z'$. If $D'$ contains $z$ and $z'$, let $D=D'$. If $D'$ contains neither $z$ nor $z'$, let $D=D'\cup \{z,z'\}$. Finally, if $D'$ contains exactly one of the centers, say $z'$, let $D=D'\cup \{z,w\}$, where $w$ is a neighbor of $z$ in $M_j$ not belonging to $D'$, that always exists; see Figure~\ref{fig:cases8_2}. In the three cases, $D$ is a semi-PD-set in $T$, semipairing $z$ and $z'$ in the second case, and $z$ and $w$ in the third one.
	
\end{enumerate}

Therefore, in all cases, we constructed a semipaired dominating set in $T$ of size at most $\cotas$, so the theorem holds.	
\end{proof}




\section*{Acknowledgments}

C. Hernando, M. Mora and J. Tejel are supported by project PID2023-150725NB-I00, funded by MICIU/AEI/10.13039/501100011033.
J. Tejel is also supported by project E41-23R,  funded by Gobierno de Arag\'on.
M. Claverol was supported by project PID2019-104129GB-I00/AEI/10.13039/501100011033 of the Spanish Ministry of Science and Innovation.


\bibliographystyle{abbrvurl}
\bibliography{bibliography}

\end{document}